\documentclass[a4paper,12pt
]{article}
\usepackage{mathpazo}
\usepackage{amsmath}
\usepackage{amssymb}
\usepackage{latexsym}
\usepackage{graphicx}
\newtheorem{thm}{Theorem}[section]
\newtheorem{lem}[thm]{Lemma}
\newtheorem{prop}[thm]{Proposition}
\newtheorem{ass}[thm]{Assertion}
\newtheorem{cor}[thm]{Corollary}
\newtheorem{rem}[thm]{Remark}
\newtheorem{df}[thm]{Definition}
\newtheorem{ex}[thm]{Example}
\newtheorem{conj}[thm]{Conjecture}
\newtheorem{prob}[thm]{Problem}
\newtheorem{fig}[thm]{Figure}
\newtheorem{statement}[thm]{}

\newcommand{\Prop}[1]{\begin{prop}{\rm\quad{#1}}\end{prop}}
\newcommand{\Ass}[1]{\begin{ass}{\rm\quad{#1}}\end{ass}}
\newcommand{\Cor}[1]{\begin{cor}{\rm\quad{#1}}\end{cor}}
\newcommand{\Rem}[1]{\begin{rem}{\rm\quad{#1}}\end{rem}}

\newcommand{\Prob}[1]{\begin{prob}{\rm\quad{#1}}\end{prob}}

\newcommand{\proof}{\noindent{\it Proof.} }
\newcommand{\proofof}[1]
{\noindent{\it Proof} of {#1}.  }
\newcommand{\QED}{Q.E.D. \hspace{\fill}$\square$}
\newcommand{\g}{{\frak g}}
\newcommand{\R}{{\Bbb R}}
\newcommand{\N}{{\Bbb N}}

\newcommand{\Z}{{\Bbb Z}}

\newcommand{\Hyp}{{\Bbb H}}

\newcommand{\D}{{\cal D}}
\newcommand{\E}{{\cal E}}

\newcommand{\F}{{\cal F}}
\newcommand{\G}{{\cal G}}

\newcommand{\IE}{{\it i.e.}, }
\newcommand{\EG}{{\it e.g.}, }
\newcommand{\RESP}[1]{[{\it resp.} {#1}]}

\newcommand{\W}{{\cal W}}

\newcommand{\ds}{\displaystyle}

\newcommand{\bul}{$\bullet$}

\newcommand{\vol}{d\mathrm{vol}}

\newlength{\itemwidth}
\newcommand{\ul}{\underline}

\begin{document}
%
\begin{center}
{\large 
{\bf \LARGE 
Geometry and dynamics 
\\
of 
\\
Engel structures
} \vspace{23pt} \\
Yoshihiko MITSUMATSU\footnote{
$\!\!\!\!\!\!$ Supported in part 
by Grant-in-Aid for 
Scientific Research (B) 17H02845
and by Chuo University Grant for Special Research
2016/2017

{\it Date}: April 17, 2018.   
 {\copyright}  2018 Y. Mitsumatsu.  


Subsections with 
${}^*$ deal with issues 
which are new at least to the author. 
}
(Chuo University, Tokyo)
}\vspace{5pt}
\end{center}
\newlength{\abst}
\setlength{\abst}{120mm}

\begin{center}
\begin{minipage}{\abst}
\indent 
 The aim of this paper is to extend 
basic understanding of Engel structures 
through developing geometric constructions 
which are canonical to a certain degree 
and 
the dynamics 
of Cauchy characteristics in the transverse spaces 
which may exhibit elliptic, parabolic, 
or hyperbolic natures in typical cases. 
\end{minipage}
\end{center}
%
%
%
%
%
%
%
%
\tableofcontents
\setcounter{section}{-1}
\section{\large Introduction}
\indent 
This article is devoted to extending 
our basic understandings of Engel structures 
from two points of view.  
As basic constructions, 
the Cartan prolongation from the contact 3-manifolds 
and the Lorentz prolongation from the Lorentzian 3-manifolds 
are known since Elie Cartan \cite{C}. 
We review these constructions 
and introduce a new one among several other methods 
of basic constructions, which is called 
the {\it pre-quantum prolongation}. 
Not only review the formality of these constructions, 
we give many examples of these constructions. 
One big class of examples are coming from 
surfaces with metrics.  
The second aim is to discuss the nature of the dynamics of the 
Cauchy characteristic acting on the even contact structure 
modulo the Cauchy characteristic. 
Kotschick and Vogel proposed the notion of 
(weak-)hyperbolicity for 
this action and raised very illustrative examples. 
We will develop slightly further this notion 
and formulate the ellipticity, 
parabolicity, or hyperbolicity of this action.  

The Cartan prolongation provides 
elliptic structures in this sense.  
Looking at hyperbolic structures or negatively curved 
metrics on surfaces, 
the Lorentz prolongation 
provides a bunch of hyperbolic ones. 
Therefore looking for ``parabolic'' structures coming 
from to a certain degree canonical constructions 
is one of the motivation and in fact 
a strong driving force of this study.  

Recently from the point of view of construction 
of Engel structures, 
after the break through by Vogel \cite{V},  
an h-principle oriented study has been developed 
.  
However the notion of Engel structure has not yet been 
widely common.  So we start with the article with some standard 
definitions and then proceed to provide basic constructions 
including classical ones, \IE the Cartan and  Lorentz prolongations, 
as well as a new one. Then in later sections after some basic 
notions on the dynamics of Engel structures have been introduced, 
more detailed observations are given to the Lorentz prolongations 
coming from surfaces with definite or non-definite Riemannian metrics 
and also to the new construction.  
In this article we assume that manifolds, distributions, 
and any other objects 
are smooth unless otherwise stated. 
\bigskip
\\
\noindent
{\bf Acknowledgement}\quad 
The author is deeply grateful 
to the American Institute of Mathematics and the organizers 
of the workshop ``Engel Structures'' held at the AIM 
from April 17 to 
21, 2017.  
It gave the author an opportunity 
that the Engel structure started to draw 
the author's strong attention. 
Also he is grateful to 
Robert Bryant from whom the author learned a lot. 
This note is a report to the workshop 
from the author. 
He is also grateful to 
to Hajime Sato for 
introducing contact geometry of 
differential equations 
and 
to 
Goo Ishikawa for knoticing 
concepts such as 
Lorentz prolongation 
and 
W\"unschmann invariants.

\section{Basic concepts and constructions 
in the study of Engel structures}
\subsection{Basic definitions}
An Engel structure $\D$ on a 4-manifold $M$ is a $2$-dimensional 
distribution which satisfies the following non-integrability condition 
on the derived distributions;  
\begin{eqnarray*}
\quad 
&(\mathrm{D1})&: \, \mathrm{rank}\,\E\equiv 3 \,\,\,\, 
\mathrm{where}\,\,\,\,  \E=[\D,\D],  \,\,\,\,
\qquad\qquad\qquad\qquad\qquad\mbox{}
\\
\mathrm{and}& & \,\,\,\, 
\qquad\qquad\qquad\qquad\qquad\mbox{}
\\
\quad 
&(\mathrm{D2})& : \, \mathrm{rank}\, [\E,\E] \equiv 4\,, 
\,\,  \IE \,\, [\E,\E]  \equiv TM. 
\qquad\qquad\qquad\qquad\qquad\mbox{}
\end{eqnarray*}
Here by abuse of notation, the distributions also denote 
the sheaves of germs of their smooth sections. 
It is well-known and easy to show that 
in general for any smooth distribution $\mathcal A$ of 
constant rank, $\mathcal A\subset[ \mathcal A,  \mathcal A]$ and 
the Lie bracket 
$[\,,\,] : \mathcal A \otimes \mathcal A \to 
 [ \mathcal A,  \mathcal A]/\mathcal A$ is 
skew symmetric and tensorial. 
Therefore the surjectivity of 
$[\,,\,] : \E \otimes \E \to  [ \E,  \E]/\E=TM/\E $ 
implies the annihilator $\W\subset \E$ 
is of constant rank one. 
Now we see that $\W$ is contained in $\D$ becuase otherwise we have 
$\D\cong \E/\W$ which implies $[\,,\,] : \D \otimes \D \to TM/\E$ 
is surjective while we have $[\D, \D]=\E$. 

This line field $\W$ is called the {\em Cauchy characteristic} 
of $\D$,  or {\em characteristic foliation},  
{\em characteristic line field}, {\it etc}... and plays a 
fundamental role in the study of Engel structures.  
Let $\F_\W$ denote the 1 dimensional foliation which 
$\W$ defines. 

Let us take a flow box neighborhood $U$ of any point $P\in M$, 
\IE 
a coordinate neighborhood $U\equiv D^4$ with local coordinate 
$(X,Y,Z,W)$ such that 
the Cauchy characteristic foliation $\F_\W$ 
coincides with the local foliation 
$\{X\equiv Y\equiv Z\equiv \mathrm{constant}\}$. 
By definition $\E$ is invariant under the flow along $\W$ and 
induces an invariant plane field on the transverse spaces 
to this foliation.  
Therefore on the local 3-dimensional quotient space 
$U/{\F_\W}=\{(x,y,z)\}$ 
a plane field $\xi$ is induced. 
This plane field is in fact a contact plane field 
because of (D2). 
For this reason, we call $\E$ the {\it even contact structure} 
associated with the Engel structure $\D$.

We can modify the local coordinate into $(x,y,z,w)$ as follows. 
First if we need we take $U$ smaller such that 
$\xi=\mathrm{ker}\,[dy-zdx]$ on $U/{\F_\W}$ 
and at the point $P$ $\D_P=\mathrm{ker}\,dy
\cap
\mathrm{ker}\,dz
$ hold, 
thanks to the conatct Darboux theorem 
(See \cite{G} for the fundamental theory for contact structures). 
If necessary we take $U$ further smaller 
and we can define the function $w$ so as to satisfy 
$\D=\mathrm{ker}\,[dz-wdx]\cap\E$.  
The condition (D1) implies 
 $dw$ never vanishes on $U$ and 
$dx$, $dy$, $dz$, and $dw$ are linearly independent.  
Therefore we can take the local coordinate $(x,y,z,w)$ 
such that $\E$ is defined by the 1-form $dy-zdx$ 
and $\D$ is defined by the pair of 1-forms
$$ 
dy -zdx\,,\qquad dz-wdx.
$$
We call this coordinate neighborhood the 
{\it Engel-Darboux} coordinate neighborhood 
(or `chart'  for short).  

As the Engel-Darboux chart suggests, 
the most important and fundamental example of 
Engel manifolds is the space $J^2(1,1)$ 
of 2-jets of smooth functions of one variable.  
As a space $J^2(1,1)$ is nothing but 
$\R^4=\{(x,y,z,w)\}$. 
Let us consider 
smooth functions $y=f(x)$ and its derivatives 
$z=f'(x)$ and $w=f''(x)$. 
Then its 2-jet graph 
$\Gamma J^2(f)=\{(x,y,z,w) = (x, f(x), f(x)', f(x)'')\}$  
is tangent to the canonical plane field $\D=\D_0$ 
defined by the pair of 1-forms 
$dy -zdx $ 
and 
$dz-wdx$.    
Let us fix a particular frame 
$$
[\mathrm{EF}]:\qquad
\ds X=\frac{\partial}{\partial x}
+z\frac{\partial}{\partial y}
+w\frac{\partial}{\partial z}, \quad
Y=\frac{\partial}{\partial y}, \quad
Z=\frac{\partial}{\partial z}, \quad
W=\frac{\partial}{\partial w}
\qquad\qquad\quad
$$
adapted to the Engel-Darboux coordinate. 
Of course the Cauchy cahracteristic of this Engel structure 
is spanned by $W$ 
and the even contact structure $\E$ is defined by 
the 1-form 
$dy -zdx$.  
Therefore the full flag associated with the Engel stucture 
is indicated as follows. 
$$
\begin{array}{ccccccccc}
0 & \!\subset\! & \!\W\! & \!\subset\! & \!\D\! & \!\subset\! 
& \!\E\! & \!\subset\! & \!T\R^4\! 
\\
 & & \Vert & &  \Vert & & \Vert & &  \Vert 
\\
 & &\langle W \rangle & & \,\,\,\W\oplus\langle X \rangle & 
& \,\,\,\,\D\oplus\langle Z \rangle & 
&\,\,\,\,\E\oplus\langle Y \rangle
\end{array}
$$
\Rem{The vector field $X$ looks just complementary, while in the 
Lorentz prolongation it has a certain importance. 
See Subsection 1.3. 
}

The coordinates $(x,y,z,\theta)$ with  
the pair of 1-forms
$$
dy -zdx \quad\mathrm{and}\quad \cos\theta\, dz- \sin\theta\,dx
$$
also define an Engel structure $D_l$ on $\R^4$.  
Its restriction  $\D_l\vert_{\R^3\times (-\pi/2, \pi/2)}$ 
is isomorphic to the standard Engel structure $\D_0$ 
through the identification $w=\tan\theta$. 
In a Engel manifold we can take local coordinates 
as above. 
Then we call them a {\it long Engel-Darboux} coordinate, 
no matter how long we can take $\theta$ along $\W$ curves. 
\medskip

Once we have an Engel structure $\D$ on $M$, 
it gives rise to a full flag 
$\W\subset \D \subset \E \subset TM$.  
We have four real line bundles none of which is 
necessarily oriented, while their orientabilities 
are related to a certain degree.  
First of all, $TM/\W$ is oriented because 
$\E/\W$ is a contact structure on this 3-dimensional space.  
$\W$ is not necessarily oriented, but if we give locally 
an orientation to $\W$, through the movement of $\D/\W$ in 
$\E/\W$ along $\W$, it also defines an orientation of 
$\E/\W$ locally. Now $TM/\E$ is spanned by the Reeb vector field 
of the contact structure $\E/\W$ on $M/\W$. 
Therefore the orientation sheaf of $TM/\E$ is also coherent 
to those of $\E/\W$ and of $\W$. 
The orientability of $\D/\W$ and hence that of $\E/\D$ are 
quite independent of others. 
Summarizing above all in other words,  
the two vector bundles of rank three, \IE 
$\E$ and $TM/\W$, 
are oriented and there is no more constraint.

\subsection{Cartan prolongation}
In this and the next subsection, 
we review two classical constructions of Engel structures 
from certain geometric structures on 3-manifolds, 
both of which Elie Cartan understood well (\cite{C}). 

The first procedure which we introduce here is called the 
{\it prolongation} or {\it Cartan prolongation}. 
Let us take a contact plane field $\xi$ on a 3-manifold $N$ 
and the projectification $\pi : M=\R P(\xi)\to N$ of $\xi$ 
which is an $\R P^1$-bundle over $N$. Instead of taking 
 $\R P^1$-bundle, also we can take the associated 
$S^1$-bundle or if $\xi$ is a topologically trivial 
$R^2$-bundle we can take its infinite cyclic covering 
which is a principal $\R^1$-bundle.  
Let  $(n,\ell)$  denote a point in $M$, namely, 
$n$ is a point in $N$ and $\ell$ is a line in $\xi_n$. 
We define the plane field $\D$ on $M$ as 
$\D_{(n,\ell)}=(D\pi)^{-1}(\ell) \subset T_{(n,\ell)}M$.

In this construction, 
it is straightforward to see from the definition 
that the even contact structure 
$\E$ is the pull-back of the contact structure 
$\xi$ on $N$ by the projection $\pi$ and 
the Cauchy characteristic $\W$ is the line field tangent to 
the fibres of $\pi$. All the characteristic lines are closed 
and the (transverse) holonomy of the 1-dimensional foliation 
$\F_\W$ along any closed leaf is trivial.

\subsection{Lorentz prolongation}\label{LorentzProlongation}
The {\it Lorentz prolongation} is reviewed here.  
Let us take a 
Lorentzian 3-manifold $(V, dg)$, \IE    
a smooth three manifold $N$ and a non-definite inner 
product on $TV$. We assume $dg$ has the signature $(2,1)$.   
A null line in $T_vV$ is a line which is null with respect to 
$dg$. It is a line on the `light cone' and thus 
the set of all such lines is a circle, 
which we call the {\it null circle} and is denoted by 
$NC(T_vV)$.  
It is worth noting that in $\R P(T_vV)$ it is a priori 
not a linear circle but a quadratic one. 
Then we consider the null-circle bundle 
 $\pi : M=NC(TV) \to V$ on which the Engel plane 
field $\D$ is defined to be 
$\D_{(v,\ell)}=(D\pi)^{-1}(\ell) \subset T_{(v,\ell)}M$.    
Here $(v,\ell)$ denotes the point in 
$NC(T_vV)$ 
indicating the null line $l\subset T_vV$.  
\begin{thm}[Lorentz prolongation,  
E. Cartan \cite{C}]\label{Lorentz Prolongation}
{\rm \quad 
1) The first derived distribution $\E=[\D,\,\D]$ is 
given as $\E_{(v,\ell)}=(D\pi)^{-1}(\ell^\perp)$. 
\\
2)  The second derived $[\E,\,\E]$ is $TM$. 
\\
3)  Therefore $\D$ is an Engel structure on 
$M=NC(TV)$. 
}
\end{thm}
The Cauchy characteristic of the Lorentz prolongation 
exhibits importance in various senses. 
\begin{thm}[
Cauchy characteristics and null-geodesics,  
E. Cartan \cite{C}]\label{null-geodesic}
{\rm \quad 
The Cauchy characteristic $\W$ of the Engel structure 
$\D$ is the line field 
given by the null-geodesic flow, namely, the natural lifts 
of null-geodesics on $V$ to $NC(TV)$. 
}
\end{thm}
\Cor{
Up to parametrization, the null-geodesics are 
invariant under conformal change of Lorentzian metrics. 
Namely, if $\gamma(t)$ is a null geodesic of 
Lore tzian manifold $(M^3, dg)$,  
then for any smoooth positive function $f$ on $M$ 
$\gamma(t(s))$ is also a null-geodesic 
for the Lorentzian metric $fdg$ 
with some reparametrization $t(s)$.  
}
\Rem{This fact is well-known even for Lorentzian manifolds 
of any dimension. In our dimension, 
we can understand it from the process of Lorentz 
prolongation, because it only depends on the 
conformal class. 
}
For the sake of being self-contained and for 
understanding constructions in this article, 
we give a proof of these theorems.  
For Theorem 1.2 another proof is given in 
\ref{NullGeodesicProof}  
which is 
based on a characterization of Cauchy characteristic curves 
in terms of infinitesimal deformation whose formulation 
is given in 
\ref{InfinitesimalRigidity}. 
\medskip

\proof 
\\
The arguments can be done locally, so take a small open set $U$ 
in $V$ and $NC(TU)$ is considered. 
First take a light-like vector field, namely a smooth vector field 
$k$ satisfying $dg(k, k)\equiv -1$. 
Next take a smooth orthonormal frame 
$\langle e$, $f\rangle$ 
of the orthogonal complement 
$\langle k \rangle^\perp$ 
which is positive definite with respect 
to the Lorentzian metric $dg$.  
Then the null-line $\ell$ is indicated as 
$\ell = \langle \cos\theta\, e + \sin\theta\, f + k\rangle$.  
So  the null circle at each point $n\in N$ 
parameterized by $\theta\in S^1$. 
Therefore $NC(TU)$ is now identified with 
$U \times S^1$ in this sense. 
By definition, 
$\D$ is spanned by two vector fields 
$F=\frac{\partial}{\partial \theta}$ 
and 
$L=\cos\theta\, e + \sin\theta\, f + k$.  
Clearly we have $Y=[F,L] =-\sin\theta\, e + \cos\theta\, f$ 
which spans the orthonormal complement 
$\langle L \rangle^\perp$ together with $L$. 
This explains 1). 
Now we have also 
$[F, Y]=- (\cos\theta\, e + \sin\theta\, f) = -X$. 
As $X$,$Y$,and $k=L-X$ span $TU$, and $T(NC(TU))$ together with $F$,  
we see that $[\E, \, \E]=TM$. This proves 2) and thus 3). 
\medskip

To prove the second theorem we need more precise computation. 
On a Lorentzian manifold, 
exactly the same in the case of Riemannian manifolds, 
there exists a unique connection $\nabla$, 
the (Lorentzian) Levi-Civita connection, which is 
compatible with the Lorentzian metric $dg$ and 
torsion free,, namely, for any vector fields $A$, $B$, $C$,
the two properties   
\begin{eqnarray*}
& & A\, dg(B,\,C)=dg(\nabla_A B,\, C) + dg(B,\,\nabla_AC), 
\\
& & \nabla_A B -\nabla_BA =[A,\, B]
\end{eqnarray*}
are satisfied. 
A curve $\gamma(t)$ on $V$ is a geodesic iff 
$\nabla_{\dot{\gamma}(t)}\dot{\gamma}(t)=0$. 
A geodesic is a null-geodesic iff 
$dg(\dot{\gamma}(t),\,\dot{\gamma}(t))=0$. 
Of course $dg(\dot{\gamma}(t),\,\dot{\gamma}(t))=0$ for some 
$t$ implies the same holds for any $t$. 
We take a small neighborhood $U$ of point in $V$ as above 
and prove the second the second theorem 
on $U\times S^1 \subset M$.  

\begin{rem}{\rm\quad 
For a null-geodesic $\gamma(t)$ with $\dot\gamma(t)\ne 0$ 
$\gamma(\varphi(t))$ is again a geodesic iff $\varphi(t)$ is 
a constant. On the other hand, it is in general impossible to 
normalize the velocity globally. In other words, 
even though a null geodesic $\gamma(t)$ comes back 
to an initial point on a closed trajectory, 
namely $\gamma(0)=\gamma(1)$ and  
$\dot\gamma(0)=c\dot\gamma(1)$ for some $c>0$, 
it does not imply  $\dot\gamma(0)=\dot\gamma(1)$. 
}
\end{rem}

Take every null-geodesics $\gamma(t)$ 
which pass through a point in $U$ and 
consider their natural lift 
$(\dot\gamma(t), \langle\dot\gamma(t)\rangle)$ to $NC(TU)$. 
By taking a smaller space-like disk $D\cong D^2$ and 
the orbits inside $U$  which meet the disk $D$, 
we can assume that such orbits fill up 
a neighborhood $\tilde{U}$ in $NC(TU)$ 
of the fibre $\pi^{-1}(v)$ of the center $v$ of $D$ 
and each orbit is a segment. 
Here $U$ is also modified to be the union of such orbits.  
Therefore we can smoothly assign the parameterization of 
each geodesics in this neighborhood $\tilde{U}$. 
Then we have a local vector field $\Gamma$ which generates 
the local null-geodesic flow on $\tilde{U}$.  
Using the local framing which we used for the above proof, 
$\tilde{U}$ can be regarded as the product $U\times S^1$. 
Then with respect to this product structure, we have 
$\Gamma(v,\theta) = \Gamma_V + f(v, \theta)F$ where 
$f$ is a smooth function.

According to this product structure, the connection 
$\tilde\nabla$ 
on $T(NC(TU))$ is also defined as the product 
$\tilde\nabla =\nabla^L \oplus \nabla^{S^1}$ 
of the 
Levi-Civita connection $\nabla^L$ on $T\tilde{U}$ 
associated with the Lorentzian metric on $U$ 
and the trivial connection $\nabla^{S^1}$ on $TS^1$. 
The connection $\tilde\nabla$ is a symmetric, \IE torsion free, 
and compatible with the product metric. It is also 
compatible with the partial metric $\pi^*dg$.  

In this formulation, the condition that $\Gamma$ generates 
a local null-geodesic flow is described that 
the  $V$-component of $\tilde{\nabla}_{\Gamma} \Gamma_V$ 
is trivial.  
This is equivalent to that 
$\tilde{\nabla}_{\Gamma} \Gamma$ has trivial $V$-component 
because $F$ is a parallel field.

The statement 4) is nothing but $[\Gamma, \E]\subset\E$, 
where $\E_{(v,\ell)}=(D\pi_{(v,\ell)})^{-1}(\ell^\perp)$. 
As we have seen that $[\Gamma, (D\pi)^{-1}(\ell)]\subset\E$,  
what we have to show is $[\Gamma, Y]\subset\E$, namely, 
$\pi^*dg(\Gamma, [\Gamma, Y])=0$,   
which is computed as follows.  
\begin{eqnarray*}
\pi^*dg(\Gamma, [\Gamma, Y])
&=& 
\pi^*dg(\Gamma, \tilde\nabla_\Gamma Y-\tilde\nabla_Y\Gamma)
\\
&=& 
\pi^*dg(\Gamma, \tilde\nabla_\Gamma Y)
-\pi^*dg(\Gamma, \tilde\nabla_Y\Gamma)
\\
&=& 
\Gamma\pi^*dg(\Gamma,  Y)
-
\pi^*dg(\tilde\nabla_\Gamma\Gamma,  Y)
-
\frac{1}{2}Y\pi^*dg(\Gamma, \Gamma)
\\
&=& 
-
\pi^*dg(\tilde\nabla_\Gamma\Gamma,  Y) =0. 
\end{eqnarray*}
\QED
\begin{rem}
{\rm\quad 1)\quad
As the second theorem is very important,  
we give an alternative proof in a later section, 
which is much shorter and relies on 
a rigidity property of Cauchy characteristic.  
It is related to the causality of the Lorentz structure. 
\\
2)\quad 
An Engel structure obtained by Lorentz prolongation is 
equipped with an extra line field 
$\langle F \rangle \subset \D$ which is transverse to 
the Cauchy characteristic $\W$ in $\D$. 
It is not true even locally that an Engel structure 
with arbitrary line field inside $\D$ which is 
transverse in $\D$ is obtained by Lorentz prolongation. 
This fact is studied by Chern in \cite{Ch} 
along an equivalence problem of 
3rd order ODE's.  
The author learned this from Robert Bryant.  
The Chern's work was initiated by W\"unschmann in his 
thesis \cite{W} under the supervision by F. Engel. 
After Chern, through the Tanaka theory, Sato-Yoshikawa 
\cite{SY} gave a geometrically clear formulation 
of the W\"unschmann invariant which is well-adapted 
to our context. 
For more historical informations, 
see \EG \cite{GN} and \cite{NP}.  

The obstruction is geometrically understood as follows. 
Locally  we can take the quotient by collapsing 
integral curves of $\langle F \rangle$ to poiints 
so that we obtain a 3-dimensional space. 
In the porjective plane of the tangent space of 
each point of this 3-dimensional space, 
the Engel plane asigns a point.  
The trace of this point along 
each $\langle F \rangle$-curve should be a 
small circle if it is obtained (locally) by 
Lorentz prolongation.  
Conversely if it is the case, 
the small circle defines a conformal class of 
Lorentzian metric on the 3-dimensional space and 
locally the Engel structure is interpreted as 
obtained from the Lorentzian prolongation. 
\\
3)\quad If we start from the flat Lorentzian space, \IE 
the Minkowski 3-space $\R^{2,1}$ the fibre direction 
$\langle F \rangle\subset \D$ coincides with 
$\langle X \rangle$ of the particular framing $[$EF$]$ 
in the Engel-Darboux coordinate 
in Subsection 1.1. 
}
\end{rem}
We will see various constructions 
by the Lorentz prolongation starting from surfaces 
with metrics.

\subsection{Pre-quantum prolongation${}^*$}
We start from richer data to construct Engel structures. 
For this, we use the construction of a 
complex line bundle with $U(1)$-connection, 
or equivalently 
an $S^1$-bundle with an $S^1$-connection 
which is well-known as the 
{\it pre-quantization} 
or 
{\it pre-quantum bundle}.  
See \EG \cite{Kos} for fundamentals. 
\begin{lem}[{
Pre-quantization}]{\rm 
For a smooth manifold $V$. 
an integral cohomology class 
$\alpha \in H^2(V; \Z)$, and 
and 
a closed 2-form $\omega$ 
which represents $\alpha$ mod torsion in 
$H^2(V;\R)$, 
there exists 
an $S^1$-bundle with an $S^1$-connection 
$\nabla$ such that its curvature 2-form $\Omega_\nabla$ 
exactly coincides with $2\pi\omega$ and the euler class   
coincides with $\alpha$. 
}
\end{lem}

Let $\xi=\mathrm{ker} \alpha$ 
be a contact structure on a 3-manifold $V$ 
with a non-singular Legendrian vector field $\underline W$ 
which is volume preserving with respect to a 
smooth volume $\vol$ and whose asymptotic cycle 
presents an integral 1st homology class. 
The last condition is also stated as 
the Poincar\'e dual closed 2-form 
$\iota_{\underline W} \vol$ 
presents an integral 2nd cohomology class in 
$H^2(V;\Z)/\!\mathrm{Torsion}\subset H^2(V;\R)$.

We take the pre-quantum $S^1$-bundle $\pi : M^4 \to V^3$ 
with a connection $\nabla$ for the closed 2-form 
$\omega=\iota_{\underline W} \vol$.   
The connection defines the horizontal hyperplane 
$H_\nabla\subset T_mM$ at each point $m\in M$.

\begin{thm}[Pre-quantum prolongation]
{\rm  1)  
The the plane field $\D$ on $M$ 
defined as 
$\D_m=H_\nabla\cap (D\pi_m)^{-1}\xi$ 
is an Engel structure. 
\\
2) The even contact structure is exactly the horizontal 
distribution $\E=H_\nabla$ and 
the Cauchy characteristic is the horizontal lift 
$\W=H_\nabla\cap (D\pi_m)^{-1}\langle {\underline W} \rangle$
of $\langle {\underline W} \rangle$. 
}
\end{thm}
\noindent
\proof
\\
The proof is divided into two parts. The first part shows 
the existence of local coordinates $(x,z,w)$ on $V$ 
which are well-adapted to 
$\xi$, $\langle\underbar{W}\rangle$, and $\omega$.  
Then in the second part 
under the preparation of the first part 
we show that the fibre coordinate $y$ can be 
so nicely chosen that $(x,y,z,w)$ forms 
an Engel-Darboux coordinate system. 
\smallskip
\\
\noindent
{\it Part I}.  
\\
In general for a non-singular vector field 
$X$ and an invariant smooth volume $\vol$ on a manifold, 
if we take a smooth positive function $f$, 
$f^{-1}\vol$ 
is also invariant under $fX$ and we have 
$\iota_{fX}f^{-1}\vol=\iota_X\vol$.   
Therefore $\iota_X\vol$, which is regarded 
as the transverse invariant volume to 
the foliation spanned by $\langle X \rangle$,  
is more essential and is closed. 

On a neighborhood $U$ of any point of $V$  
it is easy to choose smooth functions $x$ and $z$ 
on $U$ such that $dx\wedge dz=\omega|_U$ and 
that $\xi$ never coincides with $\mathrm{ker}\, dz$. 
Then on $U$, the relation $\xi= \mathrm{ker}\,[dz-wdx]$ 
defines a function $w$ on $U$. 
From the construction we see that 
$(x,z,w)$ gives a local coordinate on $U$.  
Also we have 
$\ds \left\langle\frac{\partial}{\partial w}
\right\rangle=\langle\underline{W}\rangle$, 
while 
$\ds \frac{\partial}{\partial w}$ 
does not necessarily coincides with 
$\underline{W}$ even modulo multiplication by constant. 
\smallskip
\\
{\it Part II}. 
\\
Fix a local trivialization of the pre-quantum 
$S^1$-bundle over $U$ and give the coordinate 
$(x,z,w,\theta)$ ($\theta\in S^1$). 
Then on the total space the connection 1-form 
$\Theta^\nabla$ is indicated as 
$\Theta^\nabla = d\theta + \beta(x,z,w)$ 
where $\beta(x,z,w)$ is a 1-form on the base space $V$ 
with $d\beta(x,z,w)=\omega=dx\wedge dz$. 
As $d(-zdx)=dx\wedge dz$ and by the Poincar\'e lemma, 
there exists a function $\varphi$ of $(x,z,w)$ 
such that $\beta=d\varphi -zdx$. 
Therefore we obtain 
$
\Theta^\nabla = d\theta + d\varphi -zdx 
=
d(\theta + \varphi(x,z,w)) -zdx
$.  
This implies by a gauge transformation 
$\theta \mapsto \theta + \varphi(x,z,w)$ (mod $2\pi$) 
 or 
by change of th local trivialization by 
$\theta^*=\theta + \varphi(x,z,w)$, the connection form is 
indicated as 
$\Theta^\nabla = d\theta -zdx$.  
\smallskip

For a certain part of the fibre $S^1$, 
we can give a real valued coordinate $y$ 
in place of $\theta$ 
and then we naturally obtain the followings;
\vspace{-5pt}
\begin{eqnarray*}
&\bullet& \E=H_\nabla=\mathrm{ker} [dy -zdx], 
\quad
(D\pi)^{-1}\xi =\mathrm{ker} [dz -wdx], 
\\
&\bullet& \D=H_\nabla\cap(D\pi)^{-1}\xi,
\\ 
&\bullet& \W = \left\langle\frac{\partial}{\partial w}
\right\rangle 
= (D\pi)^{-1}\langle\underbar{W}\rangle 
\cap H_\nabla\, .\vspace*{-20pt}
\end{eqnarray*}
\QED
\begin{rem}{\rm 
1) 
From a similar argument, 
the space of $S^1$-connections of a given 
$S^1$-bundle is regarded as not a vector space but 
the affine space of 1-forms. Then once the curvature form 
is specified, 
then it is the affine space of closed 1-forms. 
The difference of an exact 1-form is absorbed 
as in the above proof by a gauge transformation 
which is isotopic to the identity. 
Those which correspond to integral 1-st cohomologies 
are also absorbed by the gauge tansformations 
$M\to S^1$.  
Therefore in the above construction 
the structures might have $H^1(V;\R/2\pi\Z)$ as its moduli. 
Of course it is more complicated to consider how 
this moduli reduces 
considering diffeomorphisms of 
the total space $M$. 
\\
2) If $H^2(V;\Z)$ has torsion, the integral euler class 
is determined up to the torsion from the real class 
$[\omega]$.  Any of the $S^1$-bundles 
associated with such an integral class the pre-quantization 
works. 
Therefore provided that $V$ is compact 
we might have non-uniqueness and 
finitely many possibilities for the manifold $M$ 
from the given data $V$, $\xi$, $\vol$, and $\underline W$. 

This happens when 
we take the unit tangent bundle 
$V=S^1T\Sigma$ 
of
a closed hyperbolic surface $\Sigma$ of genus $g$. 
As $\mathrm{Tor}\,H^2(V;\Z)\cong \Z/_{2g-2}$, 
we have $(2g-2)$-many possibilities for $M$. 
\\
3) It might be worth remarking that in the pre-quantum  
prolongation, contrary to the Cartan prolongation, 
the even contact structure comes the second. 
\\
}
\end{rem}

Examples of pre-quantum prolongation are given in 
Subsection 4.4.

\subsection{Suspension by contact diffeomorphism}
\label{Suspension}

Here we review a fairy general  method 
to construct or to modify Engel structures. 

Let $(V,\xi)$ be a contact structure 
and $\ell$ be a non-singular Legendrian line field.  
Also we assume that $\xi$ is oriented plane field. 
Then we put some/any Euclidean metric on $\xi$ so that 
the angle between two Legendrian lines are defined. 
Let $\varphi$ be a contact morphism of $(V,\xi)$ 
which preserves the orientation of $\xi$.  
Also we assume that the oriented angle $d(v)$ of 
$(\varphi_*\ell)_v$ from $\ell_v$ ($v\in V)$ 
is continuously well-defined and bounded from above. 
For example if $\varphi$ is isotopic to the identity 
among contactmorphisms and $V$ is compact,     
this is the case. 
Also note that this condition is independent of 
the choice of metric. 

Let us consider the mapping torus 
$M_\varphi=\R\times V/\sim$ 
of $\varphi$, where $\sim$ identifies 
$(t+n,v)$ and $(t, \varphi^n(v))$ for $n\in\Z$. 
The contact structure $\xi$ 
it is pulled back to $\R\times V$ as a 
hyperplane field $\widetilde\E=\tilde\xi$.   
Because  $\xi$ is invariant under 
$\varphi$, 
It is also well-defined 
as a hyperplane field on $M=M_\varphi$ which 
is denoted by $\E$ 
and is going to be $\E=[\D, \D]$. 
Also let $\W$ be the suspension direction, 
\IE the natural projection image of 
$\ds \widetilde\W=\left\langle \frac{\partial}{\partial t}
\right\rangle$ 
on $\R\times V$. 

Take a smooth metric (conformal structure) on $\E/\W$ 
and pull it back to $\widetilde\E/\widetilde\W$. 
Consider the continuous twisting function $d(v)$ for 
this metric restricted to $\{0\}\times V$ and  
take an integer $K$ so that $d(v)<K\pi$ for any $v\in V$. 
On $[0,1]\times V$ let us define $\widetilde\D_{(t,v)}=
R(\rho(t,v))(\ell_v)$ where $R(\rho)$ is an rotation 
of $\widetilde\D_{(t,v)}=\xi_v$ by the angle $\rho$ 
in such a way that the smooth function $\rho(t, v)$ satisfies 
$$
\rho(0, v)\equiv 0, \quad
\rho(1, v)\equiv K\pi-d(v), \,\,\,\, 
\mathrm{and} \,\,\,\,
\frac{\partial \rho}{\partial t}>0
$$
and also that the deck transformations 
$(t+n, v)\sim (t,\varphi^n(v))$ extend 
$\widetilde\D$ on $[0,1]\times V$ to the whole 
$\widetilde\D$ on $\R\times V$ as an smooth 
Engel structure which is invariant under 
the deck transformations. 
As a result we obtain an Engel structure $\D$ 
on $M=M_\varphi$ naturally. 
\medskip

We can perform a similar modification 
to a given Engel structure $(M,\D)$. 
Let us consider a transversely embedded hypersurface 
$V\subset M$ to the Cauchy characteristic $\W$ 
with a neighborhood 
$U\cong \cup_{v\in V}[-\alpha(v), \beta(v)]\times \{v\}
\supset
\{0\}\times V=V$.  
Here, 
$\alpha$ and $\beta$ are positive smooth function on $V$,  
each curve $[-\alpha(v), \beta(v)]\times\{v\}$ 
is a Cauchy characteristic curve, and 
the first coordinate $t$ is a twisting angle 
between $(\D/\W)_{(t,v)}$ and $(\D/\W)_(0,v)$ 
with respect to some metric.   
$V$ need not be closed. 
On $V$ a contact structure $\xi=\E|_V\cup TV$ is induced. 
Consider a contactmorphism $\varphi$ of $(V,\xi)$.  
W cut the manifold $M$ along $V$ and paste it again by $\varphi$ 
to obtain a new manifold which is also denoted by $M_\varphi$. 
More precisely, remove $\{0\}\times V$, complete 
$M\setminus \{0\}\times V$ 
by the points $(0-0, v)$ and $(0+0, v)$ separately, 
and identify $(0-0, v)$ with $(0+0, \varphi(v))$ 

Because the even contact structure $\E$ is 
well-preserved by this operation, a new even contact structure 
$\E_\varphi$ is induced on $M_\varphi$. 
The new Cauchy characteristic $\W_\varphi$ is also 
naturally defined. 
Now consider the twisting function $d(v)$ for $\varphi$ 
with respect to 
$(\D/\W)_{(-\alpha(t),v)}\subset\xi_v$.   
If we can take $d(v)$ continuously so as to satisfies
$d(v)<\alpha(v)+\beta(v)$ and 
$d(v)=0$ on $V\setminus \mathrm{supp}\,\varphi$,  
there is an Engel structure on $M_\varphi$ 
which coinsids with $\D$ on $M\setminus U$ 
and whose even contact structure is $\E_\varphi$. 
It is unique up to isotopy along $\W_\varphi|_{U_\varphi}$. 

For any transversal $V$ to $\W$, 
we can find some positive functions $\alpha$ and $\beta$. 
Then contactmorphism 
which is sufficiently $C^1$-close to the identity 
with regard to $\alpha$ and $\beta$, 
the above construction is possible.  
By this modification, we can always perturb 
the dynamics of $\W$.   

\begin{ex}[Cartan prolongation]
{\rm\quad If the contact structure $\xi$ is trivial as 
an $\R^2$-bundle, it admits a non-singular Legendrian 
vector field $\ell$. Then its Cartan prolongation 
is considered to be the result of the suspension construction
by the identity, $K=1$, and $\rho(v,t)=\pi t$. 
}
\end{ex}
\begin{ex}[Bi-Engel structure, \cite{KV}]
{\rm\quad Let $\phi_t$ be a contact Anosov flow 
on a 3-manifold $V$ , namely, 
$X$ is a Reeb vector field for a contact structure $\xi$ 
and generates an Anosov flow  $\phi_t$.  
Then there is associated a bi-contact structure 
$(\xi_+, \xi_-)$, where we have 
$\xi_+\cup\xi_-=\langle X\rangle$   
and $\xi_\pm$ are twisted by $(\phi_t)_*$ in 
the opposite directions. 
Fix any $T>0$ and consider the suspension 
of $(V,\xi)$ 
by the time $T$ map $\varphi=\phi_T$ . 
Then we obtain a pair $\D_\pm$ of Engel structures 
associated with  $\xi_\pm$ for the same 
even contact structure. The twisting directions are 
opposite to each other exactly like 
the bi-contact structure. 
This example will appear repeatedly.  }
\end{ex}

\section{Action of Cauchy characteristics}
The transverse dynamics of the 1-dimensional foliation 
$\F_\W$ spanned by the Cauchy characteristic $\W$ 
is an important character of Engel structures. 

In particular, as $[\W, \E]=\E$, we can look at the action of 
$\W$ on the 2-dimensional space (vector bundle) $\E/\W$  
and the movement of $\D/\W$ inside $\E/\W$ along $\W$. 

For example, in the case of Cartan prolongation 
of a contact structure, each orbit of $\W$ is the fibre 
circle and the holonomy of $\F_\W$ is identical  
along any closed leaf, 
while 
$\D/\W$ rotates by angle $2\pi$ along each 
$\W$-orbit. 
Many more examples are studied in later sections. 

First we review the projective structure defined on each orbit 
of $\W$. This is introduced by Bryant-Hsu(\cite{BH}) and 
also studied by Inaba(\cite{I}).  

Next, we will consider to extend 
such properties of each orbit to 
the whole system $\W \curvearrowright\E/\W$.  
This is initiated by Kotschick-Vogel(\cite{KV}).

\subsection{Projective structure on Cauchy characteristic lines}
\subsubsection{Review of projective structure on 1-manifolds}
Roughly speaking, a projective structure on 1-manifold $\Lambda$ 
is a geometric structure 
modeled on $(\mathit{PGL}(2;\R), \R P^1)$ 
as a $(G,X)$-manifold, namely, 
there exists an atlas of $\Lambda$ 
whose charts take value in $\R P^1$ and the coordinate changes 
are given by elements of $\mathit{PGL}(2;\R)$. 

We do not consider non-orientable projective structure 
so that we take 
$(\mathit{PGL}^+(2;\R)=\mathit{PSL}(2;\R), \R P^1)$. 
More precisely or more formally, taking 
the universal covering 
$(\widetilde{\mathit{PSL}(2;\R)}, \widetilde{\R P^1})$ 
as the model 
and consider the developing map $\,\Phi: \tilde{\Lambda} \to 
\widetilde{\R P^1}$ which is an immersion and 
is equivariant with respect to 
$\pi_1(\Lambda)$. On $\tilde \Lambda$  
$\pi_1(\Lambda)$ acts as the covering transformation 
and on $\widetilde{\R P^1}$ through the {\it holonomy} 
homomorphism 
$\varphi:\pi_1(\Lambda) \to \widetilde{\mathit{PSL}(2;\R)}$. 
As $\Lambda$ is 1-dimensional, 
$\pi_1(\Lambda)$ is trivial or isomorphic to $\Z$. 
In the latter case, often we identify $\varphi$ with 
$\varphi(1)\in\widetilde{\mathit{PSL}(2;\R)}\,$.  
In the case of general $(G,X)$-manifolds, the developing map 
is not necessarily injective, however, so is it in our case.  
Even then it may not be surjective.   

Two projective structure on the same manifold $\Lambda$ is 
isomorphic (or called simply `the same') 
if they are united to define one atlas. 
This condition is equivalent to saying that 
if the developing maps $\Phi_1$ and  $\Phi_2$ 
are related by a single element of 
$\widetilde{\mathit{PSL}(2;\R)}$, namely, 
for some $g\in\widetilde{\mathit{PSL}(2;\R)}$ 
$\Phi_2=g\circ \Phi_1$ holds.

Two projective manifolds $\Lambda_1$ 
and  $\Lambda_2$ are isomorphic 
or diffeomorphic as projective manifolds if 
there exists a diffeomorphism 
$\psi:\Lambda_1\to \Lambda_2$ through which 
two structures are isomorphic.  
\vspace{-5pt}

\subsubsection{Projective structure on Cauchy characteristic lines
}
The projective structure of each orbit $\gamma(t)$ of $\W$ 
is defined as follows. First take the universal covering 
$\tilde\gamma$ if necessary and fix a trivialization 
$(\E/\W)|_{\tilde{\gamma}} \cong\R\times\R^2$ 
by using the action of $\W$. 
For example, using a parameterization  $\tilde\gamma(t)$ 
($t\in\R$) of the orbit, we can 
identify $(\E/\W)|_{\tilde\gamma(t)}$ with  
$(\E/\W)|_{\tilde\gamma(0)}$ 
for any $t\in \R$ by the differential 
of the holonomy transformation of $\F_\W$ along $\gamma$. 
As a foliation of codimension 3, 
a priori the linearized holonomy between 3-dimensional normal 
spaces is defined. In the case of the Cauchy characteristic 
foliation, in the normal space $TM/\W|_{\tilde\gamma(t)}$,  
$\E/\W|_{\tilde\gamma(t)}$ is 2-dimensional and 
invariant under the holonomy along $\F_\W$. 
Therefore the above identification is defined. 

Then we have the tautological developing map 
$\,\check\Phi:\tilde\gamma \to \R P^1=P((\E/\W)|_{\gamma(0)})\,$ 
defined by  
$\,\check\Phi(\tilde\gamma(t)) = (\D/\W)|_{\tilde\gamma(t)}
\in P((\E/\W)|_{\tilde\gamma(t)}\cong 
P((\E/\W)|_{\tilde\gamma(0)})$. 
The fact that $\check\Phi$ is a submersion is a direct conclusion 
of the definition of Engel structure.  
This defines a projective structure on the $\W$-orbit 
$\gamma$. 
Under this setting we do not need to take the developing map 
to $\widetilde{\R P^1}$.  
Thus the projective structure of the each orbit is 
well-defined. 

\begin{rem}{\rm  \quad 
In \cite{BH} Bryant and Hsu define the projective structure 
inside an Engel-Darboux coordinate by 
$(x,y,z,w(t))\mapsto w\in\R\cong{\R}P^1\setminus\{\infty\}$.  
Inside the Engel-Darboux coordinate a $\W$-orbit 
is exactly tangent to $\frac{\partial}{\partial w}$ 
and $x,\,y$, and $z$ are constant along it. 
Then they checked that by any Engel-Darboux chart, 
the coordinate change gives a new $w$ 
which defines the same projective structure. 

Following our definition and taking 
$\ds \frac{\partial}{\partial x}+
z\frac{\partial}{\partial y}$ 
and 
$\ds \frac{\partial}{\partial z}$ 
as the invariant frame of $\E/\W$ 
along the $\W$ orbit $(x,y,z,w(t))$ with 
$x$, $y$, and $z$ constant, 
we see that it 
coincides with the 
projective structure by Bryant-Hsu. 

In \cite{I} Inaba adopted a long 
Engel-Darboux coordinate 
$(x,y,z,\theta)$ in order to 
treat longer $\W$ orbit, which is not necessarily 
rigid. 
There, the Engel structure is defined as
$\D=\mathrm{ker}[dy-zdx] \cap 
\mathrm{ker}[\cos\theta dz - \sin\theta dx]$
and of course the developing map 
$\theta \in \widetilde{P(\E/\W)}$  
defines the projective structure. 

Remark also that inside an Engel-Darboux chart, 
we do not have to go up to the universal covering.  
}
\end{rem}
\medskip

\subsection{Projective structure on closed orbits and length}
As a typical case, here we consider 
the projective structures on closed orbits of 
the Cauchy characteristic in an Engel manifold. 

\subsubsection{Review of projective structure on circle}
\label{projective circle}
Let us review the projective structures 
on a circle $\Lambda$. 
The classification was first given by 
Kuiper (\cite{Kui}).

Consider the developing map 
between the universal coverings 
$\,\Phi:\tilde\Lambda \to \widetilde{\R P^1}\cong\R\,$ 
and its image.  
Then the holonomy 
$\,\varphi\in\widetilde{\mathit{PSL}(2;\R)}\,$ 
determines the projective structure on $\Lambda$ 
as $\,\Phi(\tilde{\Lambda})/{\varphi^\Z}$. 

To classify them, first we need to list up all the pairs 
of connected subspace of the whole line $\widetilde{\R P^1}$ 
on which the holonomy $\varphi^\Z$ acts freely. 
In an abstract sense, 
the holonomy $\varphi$ determines the structure. 
First let us recall that 
elements $A\in \mathit{PSL}(2;\R)\setminus\{\mathrm{id}\}$ 
are classified 
into three categories:
\medskip

\noindent
1)(elliptic) $|\mathrm{tr}\, A|<2$, 
no fixed point in $\R P^1$, 
conjugate to 
$\pm\begin{pmatrix}
\cos \theta & -\sin \theta \\
\sin \theta & \cos \theta
\end{pmatrix}$. 
\smallskip

\noindent
2)(parabolic) $|\mathrm{tr} \,A|=2$, 
one fixed point in $\R P^1$, 
conjugate to 
$\pm\begin{pmatrix}
1 & \mp 1 \\
0 & 1
\end{pmatrix}$.
\smallskip

\noindent
3)(hyperbolic) $|\mathrm{tr}\, A|>2$, 
two hyperbolic fixed points in $\R P^1$, 
contracting and expanding, 
conjugate to $\pm\begin{pmatrix}
e^a & 0 \\
0 & e^{-a}
\end{pmatrix}$ ($a>0$).
\medskip

Now we proceed to the classification of 
projective structures on a circle $\Lambda$. 
There are two cases, 
in one of which the image of the developing map is 
the whole line $\widetilde{\R P^1}$, while in the other case
it is an interval of finite length.  It is natural to 
regard $\R P^1=\R/\pi\Z$. The classification of 
the holonomy 
$\,\varphi\in\widetilde{\mathit{PSL}(2;\R)}\,$ 
and the projective structure on 
$\Lambda\cong S^1$ is given as follows. 
Let $\langle \varphi \rangle$ denote 
the class of $\varphi$ in $\mathit{PSL}(2;\R)$.  
\medskip

\noindent
{\bf 1) (elliptic case)}  \quad 
$|\mathrm{tr}\, \langle \varphi \rangle|<2$, 
no fixed point in $\widetilde{\R P^1}$, 
conjugate to the translation 
$t\mapsto t + (n\pi + \theta)$ 
for some $n\in \Z$, where $\theta\in (0, \pi)$ is as above. 
Basically we assume $n>0$. 
The developing image is the whole $\R=\widetilde{\R P^1}$, 
the projective length of $\Lambda$ is $n\pi + \theta$. 
The case $\langle \varphi \rangle=E$ can be also included here 
as the case of $\theta=0$, 
provided that $\varphi$ is just the translation by 
$n\pi$ for some $n\in\N$.  Apparently 
the rotations are projective symmetries, thus 
the structure is homogeneous.   
\medskip

\noindent
{\bf 2) (parabolic case)}\quad  
$|\mathrm{tr} \,\langle \varphi \rangle |=2$ 
and there are fixed points in $\widetilde{\R P^1}$. 
Then after a conjugation, it takes the following form; 
$
\langle \varphi \rangle =
\pm\begin{pmatrix}
1 & - 1 \\
0 & 1
\end{pmatrix}$, 
 the fixed point set 
Fix($\varphi$) is $\pi\Z$ and the developing image 
is $(0,\pi)$. The projective length of $\tilde\Lambda$ 
is $\pi$, while 
that of $\Lambda$ has no meaning.  
The action of $\{
\langle \varphi \rangle =
\pm\begin{pmatrix}
1 & - t \\
0 & 1
\end{pmatrix}\,\vert\,t\in\R\}$ on $\R P^1$ 
descends to 
$\Lambda$ as rotational symmetries. 
Therefore the structure is homogeneous. 
$
\pm\begin{pmatrix}
1 & + 1 \\
0 & 1
\end{pmatrix}$ 
is eliminated because it is just the inverse of the above 
and it is preferable that the action is taken 
in the positive way 
in the angle coordinate. 
\medskip

\noindent
{\bf 3) (hyperbolic case)} 
$|\mathrm{tr} \,\langle \varphi \rangle |>2$ 
and there are fixed points in $\widetilde{\R P^1}$. 
After conjugation, 
it takes the following form; 
$
\langle \varphi \rangle =
\pm\begin{pmatrix}
e^a & 0 \\
0 & e^{-a}
\end{pmatrix}$ ($a>0$), 
the fixed point set in $\widetilde{\R P^1}$ 
is Fix$(\varphi)=
\frac{1}{2}\pi\Z$, and 
the developing image is $(\pi/2, \pi)$.  
Projectively 
the length of $\tilde\Lambda$ has no meaning, because 
it can take any value in $(0,\pi)$ by conjugation.  
On the other hand, 
$|\mathrm{tr} \,\langle \varphi \rangle |>2$ 
or the derivative $\vert\log(\varphi)'\vert$ at fixed points 
(as a function $\varphi:\R\to\R$ )
is a projective invariant. 

Because the projective symmetry 
$\{\begin{pmatrix}
e^t & 0 \\
0 & e^{-t}
\end{pmatrix},\,\vert\, t\in\R\}$ 
induces the rotational symmetry on 
$\Lambda$, 
the structure is homogeneous. 
\medskip

\noindent
{\bf 4) (trans-parabolic case)}
\quad  
$|\mathrm{tr} \,\langle \varphi \rangle |=2$ 
and there are no fixed points in $\R=\widetilde{\R P^1}$. 
After conjugation let $\check\varphi\in
\widetilde{\mathit{PSL}(2;\R)}\,$ 
denote the one in the parabolic case  with 
$\langle \varphi \rangle=\langle \check\varphi \rangle
=\pm\begin{pmatrix}
1 & \mp 1 \\
0 & 1
\end{pmatrix}
$.  
Then $\varphi$ acts on $\R$ 
as $\varphi(\theta)=\check\varphi(\theta) + n \pi$ 
for some $n\in\N$, after taking the inverse if necessary.  
Thus \EG $[0, n\pi)$ is a fundamental domain of the action. 

Depending on $\mp 1$ in the off diagonal component of 
$\langle \varphi \rangle$, $n\pm$ is the projective invariant. 
The length of $\Lambda$ should be understood as 
$n\pi \pm 0$. 

On $\R$ 
the action of $\{\varphi\in
\widetilde{\mathit{PSL}(2;\R)}
\,\vert\,
\langle \varphi \rangle =
\pm\begin{pmatrix}
1 & - t \\
0 & 1
\end{pmatrix},\,t\in\R\}$ 
including those by the center 
$\theta\mapsto\theta+ k\pi\,\,k\in\Z$ 
is the symmetry. 
Therefore the decomposition 
$\Lambda=\{[\theta]\in\Lambda\,\vert\,\theta\in \pi\Z\}
\sqcup\{[\theta]\in\Lambda\,\vert\,\theta\notin \pi\Z\}
$ 
gives the orbit decomposition of the projective symmetry action. 
The structure is not homogeneous. 
\medskip

\noindent
{\bf 5) (trans-hyperbolic case)} 
\quad  
$|\mathrm{tr} \,\langle \varphi \rangle |>2$ 
and there are no fixed points in $\R=\widetilde{\R P^1}$. 
After conjugation let $\check\varphi\in
\widetilde{\mathit{PSL}(2;\R)}\,$ 
denote the one in the hyperbolic case  with 
$\langle \varphi \rangle=\langle \check\varphi \rangle
=\pm\begin{pmatrix}
e^a & 0 \\
0 & e^{-a}
\end{pmatrix}
$.  
Then $\varphi$ acts on $\R$ 
as $\varphi(\theta)=\check\varphi(\theta) + n \pi$ 
for some $n\in\N$, after taking the inverse if necessary.  
Thus \EG $[0, n\pi)$ is a fundamental domain of the action and 
$n\in \N$ and 
$\vert\mathrm{tr}\,\langle \varphi \rangle\vert$ 
or equivalently 
the derivative 
$\vert\log(\varphi)'\vert$ at $(\pi/2)\Z$ 
are the projective invariants. 
The length of $\Lambda$ should be understood as 
$n\pi \pm 0$. 

On $\R$ 
the action of $\{\varphi\in\widetilde{\mathit{PSL}(2;\R)}
\,\vert\,
\langle \varphi \rangle =
\pm\begin{pmatrix}
e^t & 0 \\
0 & e^{-t}
\end{pmatrix},\,t\in\R\}$ 
including those by the center 
$\theta\mapsto\theta+ k\pi\,\,k\in\Z$ 
is the symmetry. 
Therefore the decomposition 
$\Lambda=\{[\theta]\in\Lambda\,\vert\,\theta\in \pi\Z\}
\sqcup\{[\theta]\in\Lambda\,\vert\,\theta\in \pi(\Z+1/2)\}
\sqcup\{[\theta]\in\Lambda\,\vert\,\theta\in(0,\, \pi/2)\}
\sqcup\{[\theta]\in\Lambda\,\vert\,\theta\in (\pi/2,\,\pi)\Z\}
$ 
gives the orbit decomposition of the projective symmetry action. 
The structure is not homogeneous. 
\medskip

\subsubsection{Projective structure 
on closed Cauchy characteristic lines}
Take 
a simple closed Cauchy characteristic curve.  
Then after fixing the orientation of $\E/\W$ and $\W$ so that 
along the curve 
$\D/\W$ is moving 
in a positive angular direction. 
Then 
it naturally admits one of projective structures classified 
in \ref{projective circle}.

Take a parameterization $\Gamma:[a,b]\to M$ 
of this simple closed curve 
($\Gamma(a)=\Gamma(b)$)
a trivialization of $\E/\W|_\Gamma$ 
which is invariant under the action of $\W$ 
as mentioned in 2.1.2, in such a way that they are  
respecting the above orientations. 
The trivialization identifies $\E/\W|_{\Gamma(t)}$ with  
 $\E/\W|_{\Gamma(a)} \equiv \R^2$ for $a \leq t \leq b$. 

The differential of the first return map gives the holonomy 
$\langle \varphi \rangle \in 
\mathit{PSL}(2;\R)$. 
The (reduced) developing map 
$\langle \W/\D\rangle|_{\Gamma(t)}: [a,b]\to\R P^1=P(\E/\W|_{\Gamma(a)})$ 
lifts to the genuine developing map to $\R=\widetilde{\R P^1}$,  
which tells how many times it turns around and 
what is $\varphi \in 
\widetilde{\mathit{PSL}(2;\R)}$.  
All the ambiguities which might appear in above construction 
stay in the equivalence explained in the previous sections. 

\begin{ex}
{\rm \quad
1) (Elliptic orbit)\quad 
The Cartan prolongation $(M,\D)$ of 
a contact 3-manifold $(V,\xi)$ has 
all Cauchy characteristic closed, 
namely, they are the fibres of $\pi:M\to V$. 
Therefore the holonomy 
is trivial, it follow from the above definition 
that it is the elliptic case of length $\pi$. 

Elliptic ones with other (arbitrary) length 
naturally appear in the Lorentz 
prolongations and the pre-quantum ones as well. 
We will see them in later sections. 
\\
2) (Parabolic orbits)\quad In the Lorentz prolongation 
of the flat Lorentzian 3-torus $(T^3, dg=dx^2+dy^2-dz^2)$, 
there are lots of closed Cauchy characteristics 
corresponding to the closed geodesics of $(T^3, dg)$.  
All of them are parabolic. 
Many more example will come later. 
\\
3) Hyperbolic, trans-parabolic, or trans-hyperbolic orbits 
are easily realized by the suspension construction 
reviewed in 1.5. Of course so are elliptic or parabolic ones  
as well.   
}
\end{ex}
\begin{rem}
{\rm\quad 
In contrast with trans-hyperbolic or trans-parabolic ones, 
we may call hyperbolic or parabolic ones 
{\it genuine-hyperbolic}  or {\it genuine-parabolic}. 
In the following subsection, 
the terminology will have more meaning. 
}
\end{rem}

\bul review from Bryant-Hsu on rigidity and projective structure
and interpretation

\bul projective length

\subsection{Action of Cauchy characteristic on $\E/\W$${}^*$}

In order to understand the global structure of an Engel structure 
$\D$, 
in particular on a closed 4-manifold $M$, 
the behavior of the Cauchy characteristic $\W$ 
as 1-dimensional foliation 
and its dynamics on $M$ as well as on $\E/\W$ 
are important view points.    
Based on this the behavior of $\D/\W$ in $\E/\W$ 
along $\W$ is more clearly seen.

We reviewed in the preceding subsections that 
closed orbits of $\W$ have their own characters. 
This does not always apply to non-closed orbits 
nor to the whole structure, however, 
here we consider a very limited case where the whole structure 
still admits such a character, 
while it seems to have certain importance. 
This is similar to the differential geometric study 
of surfaces where it is not always true that the curvature 
has a single sign or is vanishing everywhere, 
while such cases have importance in various senses. 
\bigskip

Let $(M,\E)$ be an even contact structure on a closed 4-manifold 
$M$. 
We take and fix a fiberwise metric on 2-dimensional 
vector bundle $\E/\W$.  
We assume 
that $\W$ and $\E/\W$ are oriented. 
Then take a non-singular vector field $W$ which spans $\W$ 
and its flow $\phi_t=\exp tW$, whose lift to 
$\E/\W$ is denoted by $\varphi_t$. 
First we define the character of the even contact structure $\E$ 
in special cases. Remark that 
we do not have to start with an Engel structure.  
\begin{df}{\rm \quad 1) (elliptic) \quad This is the case 
where the conformal distortion of $\varphi_t$ 
is uniformly bounded.  
Precisely it is formulated as follows.  
For each point $P\in M$ 
take an oriented orthonormal basis of 
$(\E/\W)_P$, 
present the linear holonomy 
${\varphi_t}_P:(\E/\W)_P \to (\E/\W)_{\phi_t(P)}$ 
with respect to these bases, 
and let 
$\langle\varphi\rangle(t,P)$ 
denote its class in 
$\widetilde{\mathit{PGL}}^+(2;\R)
=
\widetilde{\mathit{PSL}}(2;\R)$. 
Then the uniform boundedness of the 
conformal distortion of the linear holonomy 
in $\E/\W$ is stated 
that the set 
$\{\langle\varphi\rangle(t,P)\,\vert \,
t\in \R, \, P\in M\}\subset 
\widetilde{\mathit{PSL}}(2;\R)
$ 
is bounded. If this is the case, we call 
$\E$ is of {\it elliptic} type.  
\smallskip
\\
2) (parabolic)\quad Let us assume 
there exists a real trivial sub-line bundle
$l^h$ of $\E/\W$  
and if necessary change the orientation 
of $\E/\W$.   
Then take an oriented orthonormal frame 
$\langle \ell_1, \ell_2\rangle$ 
of  $\E/\W$ 
so that  $\ell_1$ lies in $l^h$.  
If there exist positive constants $c$ and $T_0$ 
so that 
$$
\pm\frac{\langle \phi_{\pm t} {\ell_1}_P,\, 
\phi_{\pm t} {\ell_2}_P\rangle}
{\langle \phi_{\pm t} {\ell_1}_P,\, 
\phi_{\pm t} {\ell_1}_P\rangle}
\geq c\cdot t \quad \mathrm{for}\,\,
\forall P \in M, \,\, \forall t \geq T_0 
$$
the even contact structure $\E$ is said to be 
of {\it parabolic} type.  
\smallskip
\\
3) 
(hyperbolic)\quad If there exist 
two independent sub-line bundles 
$l^u$ and $l^s$ 
of $\E/\W$ which are invariant under 
the action of $\varphi_t$ and positive constants 
$c$, $c'$ and $T_0$ such that the following is satisfied.
$$
\frac
{\Vert \varphi_tv^u \Vert}
{\Vert \varphi_tv^s \Vert}
\geq c'\exp(ct)
\frac
{\Vert v^u \Vert}
{\Vert v^s \Vert}
\quad\mathrm{for}\,\,
\forall P\in M,\,\, 
\forall v^u\in l^u_P,\,\,  
\forall v^s\ne0\in l^s_P,\,\,  
\forall t\geq T_0  
$$  
Then the even contact structure is of 
{\it hyperbolic} type. 
It is also called 
{\it weakly hyperbolic} 
or sometimes {\it projectively hyperbolic}.   
Remark here that under the compactness of $M$ this 
notion is independent of the choice of fiberwise metric 
on $\E/\W$.
} 
\end{df}
\begin{rem}{\rm \quad 1)  
Note that these notions are independent of 
fiberwise metric.  
For the ellipticity (1), it is also independent of  
the choice of oriented orthonormal bases.   
\\
2)  
Moreover, in the parabolic case (2) or 
in the hyperbolic case (3), 
we can easily 
modify the fiberwise metric so that 
we can take and $T_0=0$ for (2) and (3)
and $c'=1$ for (3) 
(see \EG \cite{KV} or \cite{Mi}).  
\\
3)
Also note that in the parabolic case or in the hyperbolic case, 
there is no other $\W$-invariant (continuous) sub-line bundles 
of $\E/\W$ other than $l^h$, $l^u$, or $l^s$.   
}
\end{rem}

\begin{df}
{\rm\quad 
Let  $\D$ be an Engle structure $\D$ 
on a closed connected 4-manifold $M$.   
\smallskip
\\
(1) (elliptic) \quad 
$\D$ is of {\it elliptic} type just 
if its even contact structure $\E$ is of elliptic type,  .  
\smallskip
\\
(2) (parabolic) \quad Let us assume 
the even contact structure $\E$ 
to be parabolic.  
If  moreover $\D/\W$ does not intersect with $l^h$, 
$\D$ is called of {\it genuine-parabolic} type, 
or just of {\it parabolic} type.  
If there exists  a constant $T_1>0$ 
such that the forward orbit 
$\{\phi_t(P)\,\vert \, t\in [0, T_1]$ 
of any point $P\in M$ (\IE a portion of a $\W$-curve) 
contains a point
on which $\D/\W$ and 
$l^h$ intersect with each other , 
$\D$ is called of {\it trans-parabolic} type, 
Otherwise it is called 
 of {\it incomplete-parabolic} type.  
\smallskip
\\
(3) (hyperbolic) \quad Let us assume 
the even contact structure $\E$ 
to be hyperbolic.  
If  moreover $\D/\W$ does not intersect 
with $l^u$ nor with $l^s$, 
$\D$ is called of {\it genuine-hyperbolic} type, 
or just of {\it hyperbolic} type.  
If there exists  a constant $T_1>0$ 
such that the forward orbit 
$\{\phi_t(P)\,\vert \, t\in [0, T_1]$ 
of any point $P\in M$ (\IE a portion of a $\W$-curve) 
contains a point 
on which $\D/\W$ and 
$l^u\cup l^s$ intersect with each other, 
$\D$ is called of {\it trans-hyperbolic} type, 
Otherwise it is called 
 of {\it incomplete-hyperbolic} type.  
}
\end{df}
\begin{rem}
{\rm 
\quad
In parabolic or hyperbolic case, 
we do not know if there do exists 
incomplete ones. 
In particular, we do not know 
if there exists an Engel structure 
with parabolic even contact structure, 
which admits $\W$-orbits of two types, 
one contains a point 
on which $\D/\W$ and 
$l^u\cup l^s$ intersect with each other, 
and the other does not. 
These are fundamental problems to be considered. 

In the case of the Cartan prolongation of 
a contact 3-manifold, if the contact structure has 
trivial Euler class as plane field, 
then it is considered to be obtained 
by suspension construction 
by the identity.  Even if the Euler class 
is not trivial, 
locally  it is considered similarly and 
the resulting Engel structure is of elliptic type. 
See the following example.

If there exists a closed $\W$-curve $\Gamma$ 
in an Engel structure 
of one of the above types, 
$\Gamma$ it self has the same type. 
}
\end{rem}
\begin{prop}
{\rm\quad
Let $\E$ be a suitably oriented 
even contact structure on a closed  
4-manifold $M$ which is of parabolic or hyperbolic type.  
Then there exists an Engel structure $D$ on $M$ with 
whose even contact structure coincides with the given one 
and which is of genuine-parabolic or genuine-hyperbolic type. 
In the elliptic case, 
an even contact structure may fail 
to admit a compatible Engel structure.  
}
\end{prop}
\proof\quad In both cases, 
take the metric on $\E/\W$ as in Remark 2.5 2). 
Then it suffices to take 
$\D$ to be $\langle \ell_2\rangle\oplus\W$ in the parabolic case, 
and $\D_\pm$ to be $\langle \ell_1\pm\ell_2\rangle\oplus\W$ 
in the hyperbolic case. 
Eventually in the hyperbolic case 
we obtain a 
(genuine-hyperbolic) 
bi-Engel structure 
in the sense of Kotschick and Vogel \cite{KV}. 
\QED
\medskip

For the elliptic case, see the following examples. 
\begin{ex}{\rm \quad
Let $\xi$ be an oriented contact structure on 
$S^2\times S^1$.  Take an $S^1$-bundle 
$H=h\times\mathrm{id}_{S^1}:S^3\times S^1 \to S^2 \times S^1$ 
where $h:S^3\to S^2$ is the Hopf fibration. 
Then the even contact structure 
$\E=(DH)^{-1}\xi$ on $S^3\times S^1$ admits 
a compatible Engel structure if and only if 
$e(\xi)\ne0\in H^2(S^2\times S^1;\Z)\cong
H^2(S^2;\Z)\cong\Z$. 

In particular the standard tight contact structure on 
$S^2\times S^1$, we do not obtain an Engel structure 
in this sense.  
}
\end{ex}
\noindent
\proof
\quad
For the sake of simplicity, if $n$ is negative, 
we change the orientation of $\xi$ and assume that 
$n\geq 0$. 
If $e(\xi)=n in \Z\cong H^2(S^2\times S^1;\Z)$, 
its Cartan prolongation is defined on 
$L(2n,1)\times S^1$ whose even contact structure 
is given as $\E_{2n}=DH_{2n}^{-1}(\xi)$ 
where  $h_{2n}:L(2n,1)\to S^2$ is the $S^1$-bundle of 
Euler class $2n$ and 
$H_{2n}:L(2n,1)\times S^1 \to S^2 \times S^1$ 
is defined as 
 $H_{2n}=h_{2n}\times\mathrm{id}_{S^1}$. 
Then taking $2n$-fold covering in $L(2n,1)$ direction, 
we obtain what we want.  

On the other hand, if $n=0$, by contradiction, we 
assume that there exists an Engel structure $\D$ 
on $S^3\times S^1$ compatible with $\E$. 
Because $e(\xi)=0$, there exists non-singular 
Legendrian vector field $\ell$ on $S^2\times S^1$. 
By abuse of notation let $\ell$ also denote 
the pull-back of $\ell$ to $\E/\W$ on $S^3\times S^1$. 
Then, from the definition of Engel structures 
$\Sigma=\{P\in S^3\times S^1\,\vert\, (\D/\W)_P=\ell_P\}$ 
is a non-singular closed hypersurface of 
$S^3\times S^1$ which is 
transverse to $\W$.  
Along any fibre $H^{-1}(x)$ ($x\in S^2\times S^1$)
if we trace the movement of $\D/\W$ in $\xi_x$, 
it is clear that each fibre intersects with $\Sigma$. 
Namely $\Sigma$ is a multi-section of $H$. 
Therefore its Euler class is at most of torsion. 
This is a contradiction. 
\QED

\section{Accessible set, causality, and rigidity}
\label{AccessibleSet}
Bryant and Hsu showed that 
$W$-curves inside an Engel-Darboux coordinate neighborhood 
exhibit a rigidity property among $D$-curves. 
Inaba improved their computation and established the notion 
of accessible set, which seems perfectly fits 
into the causality property of Lorentz manifolds. 

In Subsection \ref{InfinitesimalRigidity} 
we propose an infinitesimal version of the 
rigidity 
which characterizes $W$-curves.   
This notion is valid for 
any $W$-curves of any length and 
the mechanism of the rigidity is very simple. 
Moreover it is well-adapted to give a fairy simple proof of 
Theorem \ref{null-geodesic}.  

\subsection{Rigidity of Cauchy characteristic curves and 
accessible sets}\label{Rigidity}

Bryant and Hsu found 
the following rigidity phenomena 
on the Cauchy characteristic curve 
among $\D$-curves in an Engel manifold  $(M,\D)$, 
where $\D$-curve is a smooth curve which is everywhere 
tangent to $\D$. 

Let $\gamma:[a,b]\to M$ be an embedded regular curve 
which is tangent to the Cauchy characteristic $\W$ 
and is included in an Engel-Darboux coordinate. 
Therefore after changing the coordinates of the 
Engel-Darboux chart and $[a,b]$ if necessary we may 
assume that 
$\gamma(t)=(0,0,0,t)$ for $t\in:[0, T]$ 
in an Engel-Darboux chart $\{(x,y,z,w)\}$ with 
$\D=\mathrm{ker}\,[dy-zdx]\cap
\mathrm{ker}\,[dz-wdx]$. 
We consider $\D$-curves 
which are $C^1$ close to $\gamma$.   
\begin{thm}[Bryant-Hsu, \cite{BH}]{\rm
\quad
Let $\omega(t)=(x(t), y(t), z(t), t)$ 
($0\leq t\leq T$)  
be a $\D$-curve which satisfies 
$\omega(0)=(0,0,0,0)$ and 
$\omega(T)=(0,0,0,T)$. 
Then $\omega$ coincides with $\gamma$, \IE 
$\omega(t)=(0,0,0,t)$ for $0\leq t\leq T$.  
}
\end{thm}
\proof\quad
We follow Inaba's computation (\cite{I}) 
which is easier to see. 
First let us forget the condition 
$\omega(T)=\gamma(T)$ and 
compute $y(T)$:   
$$
y(T)=\int_0^{T}\frac{dy}{dt}dt
=\int_0^{T}z\frac{dx}{dt}dt
=\int_0^{T}z\frac{1}{w}\frac{dz}{dt}dt
=\int_0^{T}\frac{1}{2w}\frac{d(z^2)}{dt}dt\, .
$$
Here the integral is ordinary, because $\ds \frac{z}{w}\to 0$ 
and $\ds\frac{1}{w}\frac{d(z^2)}{dt}\to 0$ when $t\to 0+0$.  
By integrating by parts, we have 
$$
y(T)=\frac{z^2}{2w}
-\int_0^{T}\frac{z^2}{2}\frac{d(w^{-1})}{dt}dt
=\frac{z^2}{2w}
+\int_0^{T}\frac{z^2}{2w^2}{dt}dt\,. 
$$
Therefore if we only impose $y(T)=0$ we can conclude that 
$z(t)\equiv 0$ and thus $y(t)\equiv 0$ and $x(t)\equiv 0$ 
for $t\in[0,T]$.  \QED
\\

The above computation enabled Inaba to define the 
{\it accessible set} $A$ in the Engel-Darboux chart. 

\begin{df}[Accessible set, \cite{I}]
{\rm\quad 
Let $A$ be the following subset of 
the Engel-Darboux neighborhood 
$\R^4=\{(x,y,z,w)\}$. 
\begin{eqnarray*}
A=\!\!\!\!\!\!\!\!&&A_+\cup A_- \cup A_\W
\quad \mathrm{where} \quad 
A_\W=
\{x=y=z=0\}\, ,
\\
 && A_+=
\{y>\frac{z^2}{2w}, \,\, w>0\},  \quad
A_-=
\{y<\frac{z^2}{2w}, \,\, w<0\}
.
\end{eqnarray*}
$A$ is called the {\it accessible set } 
from the origin. Note that $A_\pm$ is irrelevant 
to the $x$-coordinate. 
}
\end{df}
\begin{thm}[Inaba, \cite{I}]{\rm
\quad
1)\quad   
If a curve $\gamma:[0,T]\to\R^4$ of the form 
$\gamma(t)=(x(t), y(t), z(t), w=t)$ 
($0\leq t \leq T$) in the Engel-Darboux chart 
is a $\D$-curve which starts at the origin 
(\IE $\gamma(0)=(0,0,0,0)$), 
then the other end point $\gamma(T)$ lies 
in $A_+$ or $x(t)\equiv y(t)\equiv  z(t)\equiv 0$ 
namely $\gamma$ itself stays in $A_\W$ 
(\IE is a $\W$-curve). 

Similarly,   
a $\D$-curve $\gamma:[-T, 0]\to\R^4$ of the form 
$\gamma(t)=(x(t), y(t), z(t), w=t)$ 
($-T\leq t \leq 0$) which ends at the origin 
(\IE $\gamma(0)=(0,0,0,0)$), 
then the other end point $\gamma(-T)$ lies 
in $A_-$ or $x(t)\equiv y(t)\equiv  z(t)\equiv 0$ 
namely $\gamma$ itself stays in $A_\W$ 
(\IE is a $\W$-curve). 
\\
2)\quad
Conversely any point in $A_\pm$ can be 
joind to the origin by such a $D$-curve 
in the Engel-Darboux chart. 
If the curve touches $A_\W$, \IE 
$\gamma(t)=(0,0,0,t)$ for some $t\ne 0$, 
on $[t,0]$ or on $[0,t]$ 
(depending non the sign of $t$), 
$\gamma$ stays in $A_\W$.  
}
\end{thm}
Inaba computed the accessiboe set 
in the long Engel-Darboux chart. 
Then at the critical length $\pi$ the set is 
the natural continuation of what is described above. 
In the usual coordinate 
it is also noteworthy that the shape of the 
accessiboe set is the right cone in $(y, z, w)$-space, 
becuase 
$0=z^2-2yw=z^2+(\frac{1}{\sqrt 2}(y-w))^2 
-(\frac{1}{\sqrt 2}(z+w))^2 $.

\subsection{Infinitesimal rigidity${}^*$}
\label{InfinitesimalRigidity}
Let us introduce an infinitesimal version of the rigidity 
of Cauchy characteristic curves. 
\begin{df}\label{defLSF-IWR}
{
\rm 
For a nonsingular $\D$-curve $\gamma:[a,b]\to M$ 
is called {\it linearly strongly flexible} ({\it LSF} for short) 
iff 
there exists a smooth  deformation $\Gamma:(-\varepsilon, 
\varepsilon)\times [a,b] \to M$ for some $\varepsilon>0$ 
satisfying 
\vspace{-5pt}
\begin{eqnarray*}
(\mathrm{i})&
\Gamma(s,\cdot)=\gamma_s \,:\, 
\mathrm{non\!-\!singular}\,\, \D\mathrm{\!-curve}\,\,\, 
\forall\,s\in (-\varepsilon, \varepsilon),
\\ 
(\mathrm{ii})& 
\Gamma(s,t)=\gamma(t)
\,\,
\mathrm{on}\,\, \exists  
\mathrm{neighborhood\,\, 
of}\,\, (-\varepsilon, \varepsilon)\times\{a,b\},
\\
(\mathrm{iii})& \ds
\frac{\partial \Gamma}{\partial s}(0,t)\notin 
\E_{\Gamma(0,t)}\quad \exists t\in (a,b). 
\vspace{-10pt}
\end{eqnarray*}
$\gamma$ is called 
{\it infinitesimally weakly rigid} ({\it IWR} for short) 
if it is not LSF.   
}
\end{df}
\begin{prop}\label{propLSF-IWR}{\rm (LSF) \\
1) 
If a non-singular $\D$-curve $\gamma:[a,b]\to M$ 
is not a $\W$-curve, it is LSF. 
\\
2) 
A non-singular $\D$-curve $\gamma:[a,b]\to M$ 
is IWR iff it is a $\W$-curve, in particular, 
regardless to its length as projective structure. 
}
\end{prop}
Of course 2) implies 1) and eventually 1) can be also 
stated as ``iff''. \smallskip

\proof\quad
1) follows from the Lemma below. 
Therefore in order to prove 2) 
it is enough to show the sufficiency. 
. 

Take a $\W$-curve  $\gamma:[a,b]\to M$. 
Then we can find a long Engel-Darboux coordinates 
$(x,y,z, \theta)$ on a neighborhood of $\gamma([a,b])$ 
with 
$\D=\mathrm{ker}[dy-zdx]\cap
\mathrm{ker}[\cos\theta\,dz - \sin\theta\,dx]$
in such a way that 
$\gamma(a)=(0,0,0,0)$, $\gamma(b)=(0,0,0,\Theta)$, 
$x(t)\equiv y(t)\equiv z(t) \equiv 0$ for $t\in [a,b]$, 
are satisfied. 
So $\gamma$ is identified with 
the curve $(0,0,0,\theta
)$ for 
$\theta \in [0,\Theta]$. 
We have to show that any deformation $\Gamma$ 
with the properties (i) and (ii) 
in Definition \ref{defLSF-IWR} does not satisfies 
(iii).  
As the even contact structure $\E=\mathrm{ker}[dy-zdx]$ 
coincides with the $(x,z,\theta)$-hyperplane along $\gamma$, 
(iii) is equivalent to 
$\ds \frac{\partial y}{\partial s}(0,\theta)\ne0$. 

In order to compute 
$\ds \frac{\partial y}{\partial s}(0,\theta)$,  
we divide $[0,\Theta]$ into (possibly shorter) 
closed intervals 
so that on each interval 
$\tan \theta$ or $\cot \theta$ is well defined. 
On the intervals on which $\cot \theta$ is well-defined, 
from Inaba's computation, we have 
$$
 y(s, \theta_1) - y(s, \theta_0) 
=
\left[\frac{1}{2}z(s, \theta)^2\cot\theta
\right]_{\theta_0}^{\theta_1}
+
\int_{\theta_0}^{\theta_1}\frac{1}{2}
z(s, \theta)^2(1+\cot^2\theta)d\theta
$$
while on the intervals on which $\tan\theta$ is well-defined, 
we can exchange the roll of $x$ and $z$ by 
integration by parts and then we obtain 
\begin{eqnarray*}
 y(s, \theta_1) - y(s, \theta_0) \! 
&=&
\left[
z(s,\theta)x(s,\theta)
\right]_{\theta_0}^{\theta_1}
\\
&-&
\left[\frac{1}{2}x(s, \theta)^2\tan\theta
\right]_{\theta_0}^{\theta_1}
+
\int_{\theta_0}^{\theta_1}\frac{1}{2}
x(s, \theta)^2(1+\tan^2\theta)d\theta\,.
\end{eqnarray*}
In any case, on the right hand side each term is 
quadratic with respect to $x$ and $z$. 
Therefore if we apply 
$\ds \left.\frac{\partial}{\partial s}\right\vert_{s=0}$ 
we obtain 
$\,\,\ds
\frac{\partial y}{\partial s}(0, \theta_1)
=
\frac{\partial y}{\partial s}
(0, \theta_0) \,
$ 
because $x(0, \theta)\equiv z(0, \theta)\equiv 0$ .  
Starting from 
$\ds 
\frac{\partial y}{\partial s}(0, 0)=0$ 
and repeatedly applying the above, 
we obtain 
$\ds 
\frac{\partial y}{\partial s}(0, \theta)=0$ 
for any $\theta \in [0, \Theta]$. 
\QED
\begin{rem}{\rm 1)\quad
We see from the above argument that 
in the definition of IWR, we do not have to fix the both end 
of $\gamma$ to deform into $\Gamma$, \EG
the boundary condition 
\quad
$\Gamma(s,a)=\gamma(a)\quad   
\mathrm{for} \,\,
s\in (-\varepsilon, \varepsilon)
$
\quad 
is enough. 
\\
2)\quad
A better computation without dividing into 
intervals 
must be found 
for the proof of 2) 
if we carefully translate the computation in Lorentzian spaces 
which is done in the next subsection.  
This is left to the readers. 
}
\end{rem}

\begin{lem}[Normal form for $\D$-curves transverse to $\W$]
{\rm \mbox{}
\\
1)\quad
If a non-singular $\D$-curve $\gamma:[c,d]\to M$ 
is no where tangent to $\W$, for any $t_0\in (c,d)$ 
there is a neighborhood $[c', d']$ of $t_0$ in $[c,d]$ 
and an Engel-Darboux chart around $\gamma([c', d'])$ 
such that  $\gamma([c', d'])$ is 
an segment $[-\varepsilon, \varepsilon]$ in the $x$-axis 
for some $\varepsilon>0$ 
and $\gamma(t_0)$ is the origin. 
\\
2)\quad
Take a smooth function $f(x)$ of a single variable $x$ 
which is supported in $(-\varepsilon, \varepsilon)$ and 
satisfies $f(0)\ne 0$. Then the family of functions 
$F(s,x)=sf(x)$ for $s\in (-\varepsilon, \varepsilon)$ 
gives rise to a deformation which shows that 
$\gamma$ in 1) is LSF. 
}
\end{lem}
\proof
\quad 
The condition implies the projection of $\gamma$ to 
the local quotient space $M/\F_\W$ is an immersion. 
Take a smaller part if necessary so that it is an embedding.  
Then fix a Darboux coordinate $(x,y,z)$ for $\E/\W$ on this 
small neighborhood such that $\mathrm{ker}[dy-zdx]=\E/\W$ holds 
and the image of is as in the statement 1). 
Then naturally the statement follows. 
\QED

\subsection{Null-geodesic in Lorentzian 3-manifold${}^*$}
\label{NullGeodesicProof}
Here we give an alternative proof for 
Theorem \ref{null-geodesic} by using 
1) in Proposition \ref{propLSF-IWR}. 
\medskip

A priori, we have to show that 
a natural lift of null-geodesic is 
a $\W$-curve and the converse implication, 
while in fact it is sufficient to show only the first 
because of the following reason.  
At any point of $(v,l)\in M=NC(TV)$, 
$l$ being represented by a non-zero null vector 
$\pmb{\ell}$, 
$(\gamma(0), \dot\gamma(0))=(v, \pmb{\ell})$  
gives an initial condition of a null-geodesic 
$\gamma(t)$.  
Then the equation of (null)-geodesics admits 
at least locally a unique solution.  
It implies the natural lifts of null-geodesics 
defines a 1-dimensional smooth non-singular 
foliation on $NC(TV)$. 
\smallskip

Let $\beta:[a,b]\to V$ is a non-constant 
null-geodesic of a Lorentzian 3-manifold $(V, dg)$ 
and $(M,D)$ be its Lorentz prolongation. 
By definition, the natural lift $\gamma$ 
of $\beta$ to $M=NC(TV)$ is 
a non-singular $D$-curve. 
Take a deformation $\Gamma$  as in 
Definition \ref{defLSF-IWR} with properties 
(i) and (ii) and then consider if 
the property (iii) 
$
\frac{\partial \Gamma}{\partial s}(0,t)\notin 
\E_{\Gamma(0,t)}\quad \exists t\in (a,b)$ 
is realizable or not. 
Now the condition (iii) is equivalent to 
$$
\pi^*dg\left(\frac{\partial \Gamma}{\partial t}(0,t), 
\frac{\partial \Gamma}{\partial s}(0,t)\right) 
\ne 0 \,\,\, \mathrm{for\,\, some}\,\, t\in (a,b)\,
. $$
Here 
$\pi^*dg(\cdot,\,\cdot)$ 
denotes the pull-back of the Lorentzian metric tensor 
to $M$ by the projection $\pi:M\to V$. 
On the other hand apparently we have 
$$
\pi^*dg\left(\frac{\partial \Gamma}{\partial t}(0,t), 
\frac{\partial \Gamma}{\partial s}(0,t)\right) 
=
dg\left(\dot\beta(t), 
\frac{\partial B}{\partial s}(0,t)\right) 
$$
where $B=\pi\circ\Gamma$ 
so that we can compute them on $V$. 

For any such deformation and $t\in [a,b]$, 
using the Levi-Civita connection $\nabla$, 
we have 
\begin{eqnarray*}
&\ds dg\left(\dot\beta(t), 
\frac{\partial B}{\partial s}(0,t)\right) 
=\int_0^t
\frac{\partial }{\partial t}
dg\left(\dot\beta(\tau), 
\frac{\partial B}{\partial s}(0,\tau)\right) 
\,d\tau
\\
&\ds =
\int_0^t
dg\left(
\nabla_{\dot\beta(\tau)}
\dot\beta(\tau), 
\frac{\partial B}{\partial s}(0,\tau)\right) 
\,d\tau
+
dg\left(
\dot\beta(\tau), 
\nabla_{\dot\beta(\tau)}
\frac{\partial B}{\partial s}(0,\tau)\right) 
\,d\tau \,.
\end{eqnarray*}
The first term vanishes because $\beta$ is a geodesic. 
As we have the map 
$B:(-\varepsilon, \varepsilon)\times[a,b]\to V$ ,  
$\ds\frac{\partial B}{\partial s}$ 
and
$\ds\frac{\partial B}{\partial t}$ 
commute to each other. 
Therefore we can compute the second term as follows. 
\begin{eqnarray*}
&\ds
\int_0^t
dg\left(
\dot\beta(\tau), 
\nabla_{\dot\beta(\tau)}
\frac{\partial B}{\partial s}(0,\tau)\right) 
\,d\tau 
=
\int_0^t
dg\left(
\frac{\partial B}{\partial t}(0,\tau)
, 
\left(\nabla_{\frac{\partial B}{\partial s}}
\frac{\partial B}{\partial t}\right)(0,\tau)\right) 
d\tau
\\
&\ds
=
\int_0^t
\frac{1}{2}
\left.
\frac{\partial}{\partial s}\right\vert_{s=0}
dg\left(
\frac{\partial B}{\partial t}(s,\tau)
, 
\frac{\partial B}{\partial t}(s,\tau)\right) 
d\tau\, .
\end{eqnarray*}
The final term is nothing but  
$\,\,\ds
\frac{1}{2}
\int_0^t
\left.
\frac{\partial}{\partial s}\right\vert_{s=0}
dg\left(
\dot\beta_s(\tau), \,\dot\beta_s(\tau)
\right) 
d\tau\,\,
$ 
and vanishes because $\beta_s$ is a null-curve 
for any $s\in(-\varepsilon, \varepsilon)$.\quad  
\QED

\subsection{Causality}
We do not discuss about global causal structure but take a look at 
only local problems. 
For a general Lorentzian manifold $V$ 
and on a small neighborhood $U$ 
of a point $P\in V$, fix a future/past orientation. 
Namely, 
as the set of time-like vectors in $T_vV$ has two components 
continuous choice of one of them for each point $v$ is the 
future orientation. Then non-zero light-like(=null) 
vectors are also split into future or past oriented ones. 

Now we consider curves joining a point $P$ to another one $Q$ 
in a small neighborhood $U$ whose velocity 
is future oriented. 
It is not difficult to see that if $P$ is joined to $Q$ by 
such a curve, we can always find another curve $\omega$ 
joining $P$ to $Q$ whose velocity is always future oriented 
light-like. 

Let us fix the starting point $P$ and consider 
the local accessible set 
$A_P=\{Q\in U\, ; \, 
\gamma(0)=P,\, \gamma(1)=Q, \,\dot\gamma(t)(\leq 0):
\mathrm{future\,\,oriented}\}$ from $P$. 
It is locally a closed cone, whose interior consists 
of points to which time-like curves can reach.

the accessible set $A_P$ associated with $P$ is defined as the 
set of points at which a positive time-like curve starting from $P$ 
arrives. 
Its closure $\overline {A_P}$ is the set of points at which 
(time positive) light-like (\IE null-)curves arrive. 
If the accessible sets define a strict partial order on $V$, 
the global causality is established. 

If we consider the causality in a local sense, 
it is known that 
the boundary of accessible set is achieved 
only by null-geodesics. Together with 
the Bryant-Hsu rigidity,  this fact 
should give one more proof of 
Theorem \ref{null-geodesic}. 
This is left to the readers, while 
this motivated the proof in the previous 
subsection and the notion of IWR.

\section{Geometry and dynamics of basic examples}
In this section we look at examples 
of Lorentz prolongation and pre-quantum prolongation, 
for which the dynamical property of $\E/\W$  
is also investigated.

\subsection{Lorentz prolongation-I : Product extension${}^*$}
For a surface $\Sigma$ with a Riemannian metric $dh$, 
consider the Lorentzian $3$-manifold 
$(V,dg)=(\Sigma, dh)\times (S^1, -d\theta^2)$ 
which is just the direct product with a circle 
with negative metric. 
We call this construction the {\it product} 
({\it Lorentzian}) {\it extension}.  
Then we obtain an Engel manifold 
$(M=NC(TV), \D)$ as explained in 
\ref{LorentzProlongation}.  
Here we follow the notations there.  
We consider the the case where $\Sigma$ is complete 
and practically closed.   

The Cauchy characteristic curves are the natural 
lifts of the null-geodesics, 
while each null-geodesic of $V$ is just the combination 
of a geodesic on $\Sigma$ and that on $S^1$ 
with same speed.  
The unit tangent circle bundle $S^1(T\Sigma)$ 
admits the geodesic flow $\phi_t$, 
which is defined as 
$\phi_t((\sigma, v))=
(\gamma(t), \dot{\gamma}(t))$, 
where $\gamma$ is the unique geodesic on $\Sigma$ 
with the initial condition $\gamma(0)=\sigma$ 
and $\dot{\gamma}(0)=v$.  

Even though we do not need detailed description of 
$S^1(T\Sigma)$ and the geodesic flow, for later use, 
we fix some notations here. Let $X$ be the vector field 
which is the horizontal 
(with respect to the Levi-Civita connection) lift 
of the tautological vectors, namely, 
$\pi_*X_{(\sigma, v)}=v\in T\Sigma$. Also $Y$ denotes 
the horizontal unit vector field,  so that 
$X$ and $Y$ form an oriented orthonormal basis 
of the horizontal space. Let $Z$ denote the unit 
tangent vector field along the fibres. Therefore 
they satisfy the following commutation relations 
$$
[Z,X]=Y\, ,\quad 
[Z,Y]=-X\, ,\quad and \quad
[X,Y]= \kappa\circ\pi\cdot Z
$$
where $\kappa$ denotes the curvature function of 
$\Sigma$ and 
$\pi: S^1(T\Sigma) \to \Sigma$ is the bundle projection. 
The projection $M=NC(TV) \to V$ is also denoted by 
$\pi$ by  abuse of notation, because  
through the identification explained below 
they correspond to each other. 
The geodesic flow $\phi_t$ is generated by $X$, 
namely $X=\dot{\phi_t}\circ(\phi_t)^{-1}$ holds.  

On each point $(\sigma, \theta)\in V=\Sigma\times S^1$, 
$NC(TV_{(\sigma, \theta)})$ is identified with 
the unit tangent circle $S^1(T_\sigma\Sigma)$ 
through the identification 
$\ds v\in S^1(T_\sigma\Sigma)
\leftrightarrow
\langle
v+\frac{\partial}{\partial \theta}
\rangle
\in NC(TV_{(\sigma, \theta)})
$. 
Therefore $M=NC(TV)$ is identified with 
$S^1(T\Sigma)\times S^1$.  
Under this identification, the Cauchy characteristic is 
generated by the vector field 
$\ds X + \frac{\partial}{\partial \theta}$.   
Note that $\ds \frac{\partial}{\partial \theta}$ 
commutes with any of $X$, $Y$, and $Z$.

The first one $X^*$ from of the dual frame $X^*$,  
$Y^*$, $Z^*$ for $T^*S^1(T\Sigma)$ is, 
under the identification 
$S^1(T^*\Sigma)=S^1(T\Sigma)$ 
by the Riemannian metric $dh$,  
nothing but the tautological 1-form 
and defines the Liouville contact structure $\xi_0$, 
whose Reeb vector field is nothing but $X$. 

The Cauchy characteristic is given by 
$\ds \W=\left\langle
W=X+\frac{\partial}{\partial \theta}
\right\rangle
$
and the
Engel structure is given as the span 
$\ds \D=\left\langle
W, \, Z
\right\rangle
$.

\begin{prop}{\rm\quad
The Engel manifold $(M=NC(TV),\D)$ 
obtained as the Lorentz prolongation of 
$(\Sigma, dh)\times (S^1, -d\theta^2)$ 
is isomorphic to 
the one given by the suspension construction 
(see \ref{Suspension}) starting from 
the contact manifold 
$(S^1(T\Sigma), \xi_0=\langle Y, Z\rangle)$, 
the Legendrian field $Y$ or $Z$ 
either of which will do, and the contact diffeomorphism 
given as the time $2\pi$ map $\phi_{2\pi}$ 
of the geodesic flow, with an appropriate twisting.  
}
\end{prop}
\proof
The one by the suspension construction 
is given as follows. 
On the mapping cylinder 
$M'=S^1(T\Sigma)\times \R/\!\sim$ 
where $\sim$ is the identification 
$((\sigma, v),t+2\pi)\sim (\phi_{2\pi}((\sigma, v)), t)$, 
the Cauchy characteristic 
$\ds \W'=\left\langle
\frac{\partial}{\partial t}
\right\rangle$ 
and the even contact structure 
$\E'=\W \oplus \xi_0$ is automatically fixed.  
The Engel structure is defined as 
$\ds \D'_{((\sigma, v),t)}=\left\langle
W'\oplus  {\phi_{-t}}_*(\langle Z \rangle)
\right\rangle
$.  
On the cyclic covering 
$\tilde{M}'=S^1(T\Sigma)\times \R$, 
$\tilde{\W}'$, $\tilde{\D}'$, and $\tilde{\E}'$  
are defined as well and they are invariant under 
the deck transformation 
$T':((\sigma, v),t) \mapsto 
(\phi_{-2\pi}((\sigma, v)), t+2\pi)$.  
and thus we obtain $\D'$.  
We can check $[\D', \D']=\E'$ by the commutation 
relation $[-X, Z]=Y$ directly. But instead of doing it 
we show that $(M', \W', \D', \E')$ is isomorphic 
to the Lorentz prolongation 
 $(M', \W', \D', \E')$.   

Let us follow the identification 
$M=S^1(T\Sigma)\times S^1$ given above 
and consider its cyclic covering.  
$\tilde{M}=S^1(T\Sigma)\times \R^1$ where everything is 
lifted and indicated with $\,\tilde{\mbox{}}\,$.  
The deck transformation is 
$T:((\sigma, v), \theta) \mapsto 
((\sigma, v), \theta+2\pi)$.  
Consider the diffeomorphism 
$$
\tilde{\Phi}:\tilde{M}'\to \tilde{M},
\quad 
\tilde{\Phi}((\sigma, v), t)=(\phi_t(\sigma, v), \theta).  
$$
It is clear from the construction 
that we have 
$\tilde{\Phi}\circ T'
=T\circ \tilde{\Phi}$, 
$\tilde{\Phi}_*\tilde{\W}'=\tilde{\W}$, 
$\tilde{\Phi}_*\{\phi_{-t}*Z\}=Z$, 
and thus 
$\tilde{\Phi}_*\tilde{\D}'=\tilde{\D}$ 
as well. 
The fact that $\phi_t$ preserves 
the Liouville contact structure $\xi_0$ 
implies 
$\tilde{\Phi}_*\tilde{\E}'=\tilde{\E}$. 
Therefore $\tilde{\Phi}$ descends to the 
diffeomorphism 
$\Phi: (M', \W', \D', \E') \to (M, \W, \D, \E)$. 
\QED

Let us take a look at the dynamics of 
the Cauchy characteristic $\W$. 
$\W$ is spanned by the vector field 
$\ds W=X+ \frac{\partial}{\partial \theta}$ 
Therefore the holonomy action on $TM/\E$ is 
trivial in this construction, 
because it comes from 
$\ds \frac{\partial}{\partial \theta}$ on $S^1$. 
The action on $\E/\W$ is nothing but 
that of the geodesic flow and the curvature of the surface 
$\Sigma$ is directly reflected. 
\begin{prop}    
{\rm \quad 
In the case where the curvature $\kappa$ of $(\Sigma, dh)$ is 
positive, the action of $\W$ on $\E/\W$ is of elliptic type.   
In the case $\kappa\equiv 0$, \IE $\Sigma$ is flat, 
it is of parabolic type. 
I the case  $\kappa< 0$, it is of hyperbolic type.  
Trans-hyperbolic nor trans-parabolic case never happens 
in this construction.}
\end{prop}

\begin{rem}{\rm
\quad If we start from a surface with negative curvature, 
what we obtain is one of the bi-Engel structure of Example 1.9 
that Kotschick and Vogel obtained in \cite{KV}. 
In this case the bi-Engel structure corresponds to the 
bi-contact structure 
$(\xi_+=\langle X, Y\rangle
=\langle h, l\rangle$,  
$\xi_-=\langle X, Z\rangle
=\langle h, k\rangle)$ 
associated with the 
geodesic Anosov flow $\exp\,tX$. 
(See the next subsectioj for the notation 
$h$, $l$, and $k$.  )
In particular, 
the Engel structure we obtained here 
corresponds to $\xi_-$.  

A natural construction which gives rise to the other one 
corresponding to $\xi_+$ is given in Subsection 4.3. 
}
\end{rem}

\subsection{Lorentz prolongation-II : Magnetic extension${}^*$}
We start from a Riemannian surface $(\Sigma, dh)$ 
like in the previous subsection while  
a slightly different constrction is adopted 
to obtain a 3-dimensional 
Lorentzian manidfold $(V, dg)$. 
Let $V$ be the unit tangent circle bundle $S^1(T\Sigma)$, 
so that it admits the unique Levi-Civita connection 
$\nabla^\Sigma$. 
At each point $(\sigma, v)\!\in\! S^1(T\Sigma)$, 
the tangent space admits 
the horizontal/vertical splitting 
$\,T_{(\sigma, v)}S^1(T\Sigma) 
=H_{(\sigma, v)}\!\oplus\! V_{(\sigma, v)}$, 
where we have the natural identification 
$\,H_{(\sigma, v)}\!\!\cong\!\! T_\sigma\Sigma\,$ and 
$\,V_{(\sigma, v)}\!\!\cong\!\! T_vS^1(T_\sigma \Sigma)\,$.  
The Lorentzian metric $dg=\langle, \rangle$ 
on $V\!=\!S^1(T\Sigma)$ is 
defined as 
$\,dg\!=\!dh\!\oplus\!(-d\theta^2)\,$ 
with respect to 
the splitting. Of course $d\theta^2$ denotes 
the canonical metric of the unit circle 
$S^1(T_\sigma T\Sigma)$. 
We call this construction {\it magnetic (Lorentzian) 
extension} of a Riemannian surface $\Sigma$.  
In this subsection we are particularly interested in 
the Lorentz prolongation of the magnetic extension 
of compact srufaces with constant curvature. 

Before getting into special examples, let us fix 
some notations which are valid throughout this subsction. 
Let $X$, $Y$, and $Z$ denote the same vector fields 
on $S^1(T\Sigma)$, so that they form a Lorentzian orthonormal 
frame, namely, 
$\langle X, X \rangle
=\langle Y, Y \rangle
=-\langle Z, Z \rangle
=1$ 
and 
$\langle X, Y \rangle
=\langle Y, Z \rangle
=\langle Z, X \rangle
=0$.   
At each point $(\sigma, v)\in V=S^1(T\Sigma)$, 
the horizontal lift of $v\in T\Sigma$ to $H_{(\sigma, v)}$ 
is $X_{(\sigma, v)}$ by definition. 
$NC(T_{(\sigma, v)}V)$ is identified with 
$S^1\cong S^1(H_{(\sigma, v)})$ 
by assigning 
$X(\theta)\mapsto l=\langle X(\theta)+Z \rangle$ 
where 
$X(\theta)$ and $Y(\theta)$ 
denote 
$\cos\theta\,X + \sin \theta \, Y$ 
and 
$-\sin\theta\, X + \cos \theta\, Y$ 
respectively. 
Therefore $M=NC(TV)$ is naturally identified with 
$S^1(T\Sigma)\times S^1$.   
With respect to this product structure, 
let us introduce three horizontal vector fields 
$\tilde{X}$, $\tilde{Y}$, $\tilde{Z}$ 
and a vertical vectorfield 
$\Theta$ as 
$$
\tilde{X}\vert_{V\times\{\theta\}}=X(\theta),\,\,\,
\tilde{Y}\vert_{V\times\{\theta\}}=Y(\theta),\,\,\,
\tilde{Z}\vert_{V\times\{\theta\}}=Z,\,\,\,
\Theta=\frac{\partial}{\partial \theta}\,.
$$
which form a global frame of $TM$ 
and satisfy the commutation relations 
$$
[\tilde{Z}, \tilde{X}]=\tilde{Y},\,\, 
[\tilde{Z}, \tilde{Y}]=-\tilde{X},\,\, 
[\tilde{X}, \tilde{Y}]= \kappa\tilde{Z},\,\, 
$$
\vspace{-10pt}
$$
[\Theta, \tilde{X}]=\tilde{Y},\,\, 
[\Theta, \tilde{Y}]=-\tilde{X},\,\, 
[\Theta, \tilde{Z}]= 0 \,\,\,\,\,\,\, 
$$
and the metric relations 
$$\langle \tilde{X}, \tilde{X} \rangle
=\langle \tilde{Y}, \tilde{Y} \rangle
=-\langle \tilde{Z}, \tilde{Z} \rangle
=1,\quad 
\langle \tilde{X}, Y \rangle
=\langle \tilde{Y}, \tilde{Z} \rangle
=\langle \tilde{Z}, \tilde{X} \rangle
=0,
\vspace{-5pt}
$$
and
\vspace{-5pt}
$$
\langle\Theta,\,\cdot \rangle=0. 
$$   
Here again $\kappa$ denotes the (pull-back of) 
curvature of $\Sigma$ 
and  by abuse of notation, 
$\langle\cdot,\,\cdot \rangle$ denotes also  
the pull-back of itself 
by the projection $\pi:S^1(T\Sigma)\times S^1
\to S^1(T\Sigma)$.  
\bigskip

If the surface is a flat torus, the magnetic extension 
is the same as taking product with 
$(S^1, -d\theta^2)$.  
If the curvature $\kappa$ is not identically zero, a priori 
the resullt is different from the product extension. 

Particularly interesting is the case of hyperbolic surfaces, 
\IE $\kappa\equiv -1$. 
Let us first look at this case because this case can be 
described in a special and totally different way and 
also because even $\kappa$ is negative constant 
only the case $\kappa\equiv -1$ exhibits 
quite a different feature.

\begin{ex}[Magnetic Lorentzian extension of 
hyperbolic surface]{\rm
\quad
A hyperbolic surface is a quotient 
$\Gamma\!\setminus\!\Hyp^2$ 
where 
$\pi_1(\Sigma)\cong \Gamma\subset \mathrm{Isom}^+(\Hyp^2)
=\mathit{PSL}(2;\R)
$
is a torsion free co-compact discrete subgroup.  
The hyperbolic plane $\Hyp^2$ is described as 
$\Hyp^2=\mathit{PSL}(2;\R)/\mathit{PSO}(2;\R)$ 
and the unit tangent bundles are described as  
$S^1(T\Hyp^2)=\mathit{PSL}(2;\R)$  and 
$S^1(T\Sigma)=\Gamma\!\setminus\!\mathit{PSL}(2;\R)$.   
In its Lie algebra $\mathit{psl}(2;\R)$, 
take the basis  $h\!=${\scriptsize $\ds 
\,\frac{1}{2}
\begin{pmatrix}
1 & 0 \\ 
0 & -1
\end{pmatrix}
$}, 
$l\!=${\footnotesize 
$\ds \,\frac{1}{2}
\begin{pmatrix}
0 & 1 \\ 
1 & 0
\end{pmatrix}
$}, 
and 
$k\!=$
{\scriptsize 
$\ds \,
\frac{1}{2}
\begin{pmatrix}
0 & -1 \\ 
1 & 0
\end{pmatrix}
$}, so that each of them generates 
the 1-parameter subgroups 
{\scriptsize 
$\ds \left\{
\begin{pmatrix}
e^{t/2} & 0 \\ 
0 & e^{-t/2}
\end{pmatrix}
\right\}
$}, 
{\scriptsize 
$\ds \left\{
\begin{pmatrix}
\cosh t  & \sinh t \\ 
\sinh t & \cosh t
\end{pmatrix}
\right\}
$}, 
and 
{\scriptsize 
$\ds \left\{
\begin{pmatrix}
\cos t/2  & -\sin t/2 \\ 
\sin t/2 & \cos t/2
\end{pmatrix}{\Huge{/\{\pm 1\}}}
\right\}
$} respectively.  
As the elemets of  $\mathit{psl}(2;\R)$ 
are considered to be left-invariant 
vector fields on  $\mathit{PSL}(2;\R)$, 
they descend to 
$S^1(T\Sigma)=\Gamma\!\setminus\!\mathit{PSL}(2;\R)$. 
In this context $h$, $l$, and $k$ correspond to 
$X$, $Y$, and $Z$ respectively.

Then the left invariant vector field 
$L=X+Z \in \mathit{psl}(2;\R)$ 
canonically assigns a null-vector 
to any point $(\sigma, v) \in S^1(T\Sigma)$. 
}
\end{ex}

\begin{prop}\label{AnosovElement}
{\rm \quad
The null-vector field $L$ fills up $S^1(T\Sigma)$ 
with the null-geodesic orbits. 
}
\end{prop}
This fact is understood in many ways. 
For example, the Lorentzian metric 
on $\mathit{PSL}(2;\R)$ in this case 
is not only left-invariant 
but also right-invariant. 
Then any 1-parameter subgroup 
is a geodesic even respecting the parameter.   
Then any left translation of it is also a geodesic. 

Verifying ${\nabla}_L L =0$ by computing  
$\langle {\nabla}_L L, \cdot \rangle=0$ 
for a global framing, \EG $X$, $Y$,  and  $L$ is 
another way. A similar comlutation in more general 
setting will be done in the next example.   
A computation of this type also proves that 
a 1-parameter subgroup is a geodesic 
for a bi-invariant metric. 
\medskip

Proposition \ref{AnosovElement} implies that 
the natural lift of 
the null-geodesic which is an orbit of $L$ on 
$S^1(T\Sigma)$ lies exactly on 
$S^1(T\Sigma)\times \{\theta=0\}$. 
Here remark that if we regard 
$Y$ generating a geodesic flow of $\Sigma$ insead of $X$, 
$Y$ is of Anosov and $\langle L \rangle$ is eactly 
its stable foliation, because $[Y, L]=L$.

For $\theta\in S^1$ put 
$L(\theta)=X(\theta) + Z$ and extend $L$ to 
$\tilde L = \{L(\theta)\}_\theta$ on 
$M=NC(T(S^1(T\Sigma)))=S^1(T\Sigma)\times S^1$. 
$L(\theta)$ is the right-translation 
(or the Adjoint image) of $L$ by $\exp(\theta Z)$. 
Thus we obtain all null-geodesics in this way. 
An orbit of $L(\theta)$ lies on 
$S^1(T\Sigma)\times \{\theta\}$ and 
the vector field $W=\tilde L$ generates the Cauchy 
characteristic $\W$.

\begin{prop}{\rm \quad
For the Lorentz prolongation $\D$ 
of the magnetic extension of 
a hyperbolic surface $\Sigma$, 
the Cauchy 
characteristic $\W$ is regarded as 
an $S^1$-family of the Anosov strong stable foliations 
associated with the geodesic flow of $\Sigma$.    

In particular, it is of genuine-parabolic type.  
The invariant sub-line bundle in $\E/\W$ 
is generated by $\tilde Y$, 
the $S^1$-family of geodesic flows. 
}
\end{prop}
\proof \quad We verify the second statement. 
The commutation relation  
$[Y(\theta), L(\theta)]=L(\theta)$ 
on each $S^1(T\Sigma)\times \{\theta\}$ 
implies that 
the plane fied spanned by $Y(\theta)$ and $L(\theta)$ 
is integrable and in fact it is nothing but 
the Anosov stable foliation 
$\F^s(\theta)$ 
of the geodesic flow generated by $Y(\theta)$.  
We have seen 
$\W=\langle W \rangle=\langle {\tilde L} \rangle$,  
$\D=\langle {\tilde L}, \Theta \rangle$, 
and $\E=\langle {\tilde L}, \Theta,{\tilde Y} \rangle$.

The action of $W=\tilde L$ on 
$\E/\W$ is easily computed as 
$$
[{\tilde L}, \Theta] = -{\tilde Y}, \quad 
[{\tilde L},{\tilde Y}]={\tilde L}\equiv 0 
\,\,\,\mathrm{in} \,\,\, \E/\W.  
$$
The first one is a general phnomenon, while the second 
one is characteristic in this case. 
Along $W$, $ \theta$ inclines 
towards $\pm \tilde Y$ but it never reaches. 
This proves the $\langle {\tilde Y} \rangle$ 
is an invariant sub-line bundle of $\E/\W$ 
with which $\D/\W=\langle\Theta\rangle$  
never coincides. 
\QED 

\begin{rem}{\rm \quad 
The 2-dimensional foliation ${\tilde\F}=\{\F^s(\theta)\}_\theta$ 
on 
$S^1(T\Sigma)\times S^1$ 
naturally extends to two different 3-dimensional foliations 
$\mathcal G_0=\{S^1(T\Sigma)\times\{\theta\}$ and 
$\mathcal G_1$, whose intersection 
is exactly $\tilde \F$. 
Because 
$\F^s(\theta)=\exp(\theta Z)^*\F^s=\exp(-\theta Z)_*\F^s$, 
the leaves of $G_1$ is the trace of a leaf of $\F^s$ 
by the $S^1$ action generated by $-{\tilde Z} + \Theta$.  
Indeed, apparently we have 
$[-{\tilde Z} + \Theta, {\tilde Y}]
=[-{\tilde Z} + \Theta, {\tilde L}]=0$.  
} 
\end{rem}

\begin{ex}[Lorentzian extension of flat torus]{\rm\quad
If the surface $\Sigma$ is a flat torus, namely 
the case of $\kappa\equiv 0$, 
the magnetic extension and the product extension coincide 
to each other. In this case also the Engel structure 
by Lorentz prolongation is 
of genuine-parabolic type. } 
\end{ex}

Now we proceed to more general case. 
The goal of this subsection is the following result. 

\begin{thm}{\rm\quad 
Let $(\Sigma, dh)$ is a compact (or complete) 
Riemannian surface with constant curvature $\kappa$. 
The Engel structure by Lorentz prolongation 
of the magnetic extension of $(\Sigma, dh)$ 
is of 
\vspace{-2pt}
\begin{enumerate}
\item[(1)]\quad 
elliptic type iff $\kappa>0$ or $\kappa<-1$
\vspace{-4pt}
\item[(2)]\quad
 genuine-parabolic type iff $\kappa=0$ or $-1$
\vspace{-4pt}
\item[(3)]\quad 
genuine-hyperbolic type iff $-1<\kappa<0$. 
\vspace{-2pt}
\end{enumerate}
The sign of the quadratic function $\kappa(\kappa +1)$ 
controlls this phenomemon. 
}
\end{thm}

Apart from Engel structures, as a problem of 
magnetic Lorentzian extensions of Riemaniann surface in general, 
the following results are of certain interest.  
Also it is fundamental to understand the above theorem. 

\begin{prop}\label{NullRiemannianGeodesics}{\rm\quad 
Let $(\Sigma, dh)$ be any Riemannian surface 
(the curvature $\kappa$ can vary) 
and $(V=S^1(T\Sigma), dg)$ is its magnetic Lorentzian extension. 
\vspace{-2pt}
\begin{enumerate}
\item[1)]\quad Any null-geodesic $\Gamma(t)$ 
is of constant speed if projected down to $\Sigma$. 
\vspace{-4pt}
\item[2)]\quad Let $\gamma(t)$ be the projected image of a null-geodesic 
$\Gamma(t)$ as in 1). Then $\gamma$ is a curve with 
geodesic curvature $-\kappa(\gamma(t))$.  
\end{enumerate}
}
\end{prop}

One consequence of the above proposition 
and its proof in Engel structure 
is the following. 

\Cor{
If the curvature $\kappa$ is constant, the natural lift of 
any null-geodesic stays in a single 
$S^1(T\Sigma)\times\{\theta\}$. 
}

\proofof{Proposition 
\ref{NullRiemannianGeodesics}}
\quad
Let $(\Sigma, dh)$ be any Riemannian surface 
and $\kappa$ its (Gaussian) curvature.   
The unique Levi-Civita connection of 
the magnetic Lorentzian extension 
$(V=S^1(T\Sigma), dg)$ 
is denoted by $\nabla$. 

Though the statements are described on $V$, 
it is easier to prove them on $M=NC(TV)$. 
Therefore on $M=V\times S^1$, 
we extend $\nabla$ as 
$\nabla \times \Theta$ 
and by abuseof notation $\nabla$ denotes it again. 
Also we pull back $dg$ to $M$ 
and $\langle\, , \, \rangle$ denotes both  $dg$ on 
$V$ and the pulled-back on $M$ 
as in the proof of Proposition 4.6. 

Any point $(v, l) \in NC(TV)$ provides 
an initial condition for a null-geodesic $\Gamma(t)$ 
as $\Gamma(0)=v$, $\dot\Gamma(0)=X(\theta(0))+Z\in l$ 
which admits a unique solution $\Gamma(t)$. 
As it is already explained 
in the second paragraph of Subsection \ref{NullGeodesicProof}, 
the unique existence for the initial condition implies 
$NC(TV)$ admits a line field $\W$, which is nothing but the 
Cauchy characteristic of the Engel structure, 
whose integral curves are the natural lifts of null-geodesics. 
A priori, we do not have natural choice of vector field 
spanning $\W$. 
Therefore we consider the problem locally. 
Take a point $(v_*,l_*)\in M$ and 
a local transversal $T\cong \mathrm{int}\,D^3$ to $\W$ 
which contains  $(v_*,l_*)$.  
On $T$ the initial conditions for the 
geodesics is given as $\Gamma(0)=v\in T$, 
$\dot\Gamma(0)=X(\theta(0))+Z\in \W_v$.   
Then, for a small $\varepsilon>0$ 
$U=\{(\Gamma(t),\dot\Gamma(t))\, ; \, 
(\Gamma(0),\dot\Gamma(0))\in T,\,
\vert t \vert < \varepsilon
\}$ is an open set of $M$ 
which contains  $(v_*,l_*)$ 
and is diffeomorphic to 
$T\times (-\varepsilon, \varepsilon)$ 
and the natural lifts 
of geodesics with such initial conditions 
define a local flow on $U$, 
whose velocity field is denoted by $W$. 
On $U$, $W$ takes the following form. 
\vspace{-3pt}
$$
W =r ({\tilde X}+{\tilde Z}) 
+ f\Theta
\vspace{-3pt}
$$
Here $r$ is a positive smooth function on $U$ 
satisfying $r\vert_T\equiv 1$ and $f$ is a smooth function.  

%

The geodesic equation 
$\nabla_{\dot\Gamma}\dot\Gamma=0$ 
for a null-geodesic  $\Gamma(t)$ on $V$ 
lifts to $M$ and  $W$ satisfies 
$\nabla_W W=0$ on $U$. 
Therefore it is equivalent to 
$\langle  \nabla_W W , \, F\rangle$ 
for each member $F$ of a framing of $\pi^*V$ on $U$, 
\EG for 
${\tilde X}$,  ${\tilde Y}$,  
and  
${\tilde L}={\tilde X}+{\tilde Z}$.  
For 1) $F={\tilde X}$ safices. 
We have 
\vspace{-3pt}
\begin{eqnarray*}
0\!&=&\!\langle  \nabla_W W , \, {\tilde X} \rangle
=W\langle W ,\, {\tilde X}   \rangle
-\langle W, \, \nabla_W {\tilde X} \rangle
\\
&=&\!W\cdot r 
 -\langle W, \, [W,  {\tilde X}]\rangle
 -\langle W, \, \nabla_{\tilde X} W \rangle
\\
&=&\!
W\cdot r 
 -\langle W, \, [W,  {\tilde X}]\rangle
 -\frac{1}{2}{\tilde X}
\langle W, \, W \rangle
=
W\cdot r 
 -\langle W, \, [W,  {\tilde X}]\rangle
\vspace{-3pt}
\end{eqnarray*}
while 
\vspace{-3pt}
\begin{eqnarray*}
[W, {\tilde X}]
\!&=&\![r{\tilde X}+r{\tilde Z}+ 
f\Theta,\,
{\tilde X}]
=[r{\tilde X}, \, {\tilde X}]
+
[r{\tilde Z},\,  {\tilde X}]
+ 
[fTheta,\,
{\tilde X}]
\\
\!&=&\!-({\tilde X}\cdot r) {\tilde X} 
+
(r+f){\tilde Y}
-
({\tilde X}\cdot r) {\tilde Z}
+
({\tilde X}\cdot f) \Theta
\vspace{-3pt}
\end{eqnarray*}
implies 
$
\langle W, \, [W,  {\tilde X}]\rangle=0 \, , 
$
so that we can conlude $W\cdot r=0$, namely, 
we can take $r\equiv 1$ on $U$ and eventually on $M$. 
This completes the proof of 1). 
\medskip

Let us proceed to prove 2), which is done by 
a similar computation for $F={\tilde Y}$.  
Form 1) we can assume that globally on $M$ 
\vspace{-3pt}
$$
W ={\tilde X}+{\tilde Z} 
+ f\Theta
\vspace{-3pt}
$$
generates the null-geodesic flow. 
\Ass{
$f=-(1+\kappa)$. 
}
\proofof{Assertion} \quad We have 
\vspace{-3pt}
\begin{eqnarray*}
0\!&=&\!\langle  \nabla_W W , \, {\tilde Y} \rangle
=W\langle W ,\, {\tilde Y}   \rangle
-\langle W, \, \nabla_W {\tilde Y} \rangle
\\
&=&\!
 -\langle W, \, [W,  {\tilde Y}]\rangle
 -\langle W, \, \nabla_{\tilde Y} W \rangle
=
 -\langle W, \, [W,  {\tilde Y}]\rangle
\vspace{-5pt}
\end{eqnarray*}
and 
\vspace{-3pt}
\begin{eqnarray*}
[W, {\tilde Y}]
=
[{\tilde X}+{\tilde Z}+ 
f\Theta,\,
{\tilde Y}]
\!&=&\!
[{\tilde X}, \, {\tilde Y}]
+
[{\tilde Z},\,  {\tilde Y}]
-
f {\tilde X}
-
({\tilde Y}\cdot f) \Theta
\\
\!&=&\!
\kappa {\tilde Z} 
-
(1+f){\tilde X}
-
({\tilde Y}\cdot f) \Theta . 
\vspace{-3pt} 
\end{eqnarray*}
Therefore we obtain 
\vspace{-3pt}
\begin{eqnarray*}
0
=
-\langle W, \, [W,  {\tilde Y}]\rangle
\!&=&\!
-\langle 
{\tilde X}+{\tilde Z} 
+ f\Theta
, \, 
\kappa {\tilde Z} 
-
(1+f){\tilde X}
-
({\tilde Y}\cdot f) \Theta 
\rangle
\\
\!&=&\!
\kappa + 1 + f \, .
\vspace{-3pt} 
\end{eqnarray*}
\QED
\medskip Assertion 4.12.

Now the statement 2) follows from the following 
Proposition. 
\\
\QED Prpopsition 4.10. 
\Prop{
Let $\gamma(t)$ be a curve on $\Sigma$ 
with unit speed,   
$\Gamma(t)$ be a null-curve lift of $\gamma$ 
to $V=S^1(T\Sigma)$, 
and
$\tilde\Gamma(t)$ be the natural lift 
of $\Gamma$ to $M=NC(TV)$, 
namely, 
\medskip
\\
\qquad \bul 
$\,\,p\circ \Gamma=\gamma$, 
$\,\,\,\pi\circ \tilde\Gamma=\Gamma$ 
\\
\quad
where 
$p$ and $\pi$ denote the projections 
$V\to \Sigma$ 
and $M \to V$ respectively, 
\smallskip
\\
\qquad \bul $\,\,\Vert \dot\gamma(t) \Vert\equiv 1$, 
$\,\,\,\dot\Gamma(t)=$ a horizontal lift of 
$\dot\gamma(t) + Z$, 
$\,\,\,\tilde\Gamma(t)
=\langle \dot\Gamma(t) \rangle$. 
\medskip
\\
We present the velocity of 
$\tilde\Gamma$ is presented as 
$\dot{\tilde\Gamma}(t)
= \dot{\Gamma}(t) + 
  \varphi(t)\Theta$,  
respecting the product structure 
$M=V\times S^1$.  
Then the geodesic curvature of $\gamma(t)$ on 
$\Sigma$ is equal to $\varphi(t)+1$. 
In the integral form, 
if $\tilde\Gamma(t)= (\Gamma(t), \Phi(t))$, 
then the geodesic curvature is $\dot\Phi(t)+1$.  
}
\Rem{
1)\quad 
Proposition 4.10, in particular the statement 1) 
implies somehow the Magnetic Lorentzian extension 
remembers the Riemannian metric of $\Sigma$. 
If we look at the global symmetry 
of the magnetic Lorentzian extension, , 
unless the surface is very special type, 
\EG 
the global isometry group is the standard $S^1$-action 
in the fibre direction of the projection 
$p:S^1(T\Sigma) \to \Sigma$.  

But even as a local geometry, 
we can find 
a reminiscence of the surface 
as in the next assertion, 
which we will partly use in the proof of the above proposition.  
\\
2)\quad 
It is seen from the proof that 
the correction term $+1$ in the above proposition 
is the $+Z$ 
which makes the null-lift going up in $Z$-direction.  
}

In the magnetic extension, if we take the sum of 
the metrics of the horizontal and vertical spaces 
$H_{(\sigma, v)}$ and $V_{(\sigma, v)}$, 
both with the positive sign, 
we obtain the standard  Riemannian metric on $S^1(T\Sigma)$, 
\IE the magnetic Riemannian extension of $\Sigma$.  
\Ass{1)\quad 
Any horizontal lift of a geodesic on $\Sigma$ is again 
a geodesic on $S^1(T\Sigma)$ with magnetic Riemannian 
extension. 
\\
2)\quad Consequently the horizontal lift of 
any geodesic on $\Sigma$ is also a geodesic 
for the magnetic Lorentzian extension.  
\\
3)\quad 
The same applies to the fibre circles, 
\IE the $Z$-curves.  
}
\proofof{Assertion 4.15}
1) is easily understood because the horizontal lift is 
apparently locally minimizing the length.  
Then 2) is concluded from the fact 
that the infinitesimal 
deformation of the given horizontal lift 
and the derivation of the energy 
is computed by decomposing the deformation 
in the horizontal and the vertical directions. 
\QED{Assertion 4.15.}
\bigskip
\\
\proofof{Proposition 4.13}
First let us verify the proposition in the case of a geodesic 
$\gamma$ on $\Sigma$ of unit speed. 
Let $\zeta(\theta)$ denote the standard $S^1$-action 
in the fibre direction of $S^1(T\Sigma)$.  
Let $\Gamma_h(t)$ denote a horizontal lift of $\gamma(t)$. 
Its velocity $\dot\Gamma_h(t)$ is described as 
$X(\theta(t)_{\Gamma_h(t)})$. 
Then it is clear that $\theta_{\Gamma_h}=\theta(t)$ 
is constant in $t$. 
Any null-lift $\Gamma(t)$ of $\gamma(t)$ is 
then given as $\Gamma(t)=\zeta(t-c)(\Gamma_h(t))$ 
for some constant $c$.  
Then clearly the velocity satisfies 
 $\dot\Gamma(t)=X(\theta_\Gamma-t+c)_{\Gamma(t)}$.  
Its natural lift to $M=NC(TV)=V\times S^1$ 
is $(\Gamma_n(t), \, \Phi(t)=\theta_\Gamma-t+c)$  
and thus $\dot\Phi(t)\equiv -1$ for any geodesic 
$\gamma(t)$ on $\Sigma$. 
\smallskip

Now we consider the general case.  
The geodesic curvature $\kappa_g(t)$ 
of $\gamma(t)$ 
is defined to be 
$\nabla^\Sigma_{\dot\gamma(t)}\dot\gamma(t)$ 
or numerically to be 
$dh( \nabla^\Sigma_{\dot\gamma(t)}\dot\gamma(t), 
\nu(t))$.  
Here, $\nu(t)$ denotes the unit normal vector field
 along  $\gamma(t)$ such that 
$\dot\gamma(t)$ and $\dot\nu(t)$ form an oriented orthonormal 
basis of $T_\gamma(t)\Sigma$ 
with the Riemannian metric $dh$, and its Levi-Civita connection 
$\nabla^\Sigma$.  

Let $\Gamma_h(t)$ and $\Gamma_n(t)$ be horizontal and null- lifts 
of $\gamma(t)$ respectively.  
Present the velocity of $\Gamma_h(t)$ and of $\Gamma_n(t)$ 
as 
$\dot\Gamma_h(t)=X(\theta_h(t))_{\Gamma_h(t)}$. 
$\dot\Gamma_n(t)=X(\theta_n(t))_{\Gamma_n(t)}$. 
Then for some constant $c$, we have 
$\Gamma_n(t)=\zeta(t-c)\Gamma_h(t)$ and 
thus 
$\theta_n(t)=\theta_h(t)-t+c$. 
Therefore it is enough to show that 
$$
\dot\theta_h(t)=\kappa_g(t).  
$$

We prove the above equality at $t=t_0$.  
The following computations and arguments are done locally. 
Take the vector field $X(\theta_h(t_0))$ and consider 
an orbit $\Omega(t)$ with $\Omega(t_0)=\Gamma_h(t_0)$, 
so that $\dot\Omega(t_0)=\dot\Gamma_h(t_0)$. 
Take a local triviality of the $S^1$-bundle 
$S^1(T\Sigma) \to \Sigma$ 
so that we have a product neighborhood 
$U\times S^1$ of $\Gamma_h(t_0)$, 
$\gamma(t_0)\in U \subset \Sigma$, 
$\Gamma_h(t_0) = (\gamma(t_0), \theta_0)$.  
Also we can assume that all orbits of $X(\theta_h(t_0))$ 
in this product neighborhoos is horizontal 
with respect to this product structure. 
Then we also take a product connection 
$\nabla^P = \nabla^\Sigma\times 
(d\theta\otimes\frac{d}{d\theta})$ 
of $T(U\times S^1)$ with 
respect to this product structure. 
Then it is clear from the definition that 
$$
\langle \nabla^P_{\dot\Gamma_h(t)}\dot\Gamma_h(t), 
Y(\theta(t))\rangle = \kappa_g(t)
$$
holds. 

Now we calim that even if we replace $\nabla^P$ in this 
equation with the Levi-Civita conection 
$\nabla$ of the magnetic Lorentzian extension, 
the following argument shows that the equation still holds. 
(Instead of Lorentzian extension, we can use the magnetic 
Riemannian extension and its Livi-Civita connection. 
Even then exactly the same argument holds.) 

From the construction we have 
$
\nabla^P_{\dot\Omega(t)}\dot\Omega(t)=
\nabla_{\dot\Omega(t)}\dot\Omega(t)=0
$.  
Also it is clear that both  
$\nabla^P$ and $\nabla$ 
are torsion free as affine connection. 
Therefore for vector fiedls $A$ and $B$ 
with $B_Q=0$ at a point $Q$, 
$
(\nabla^P_AB)_Q = 
(\nabla_AB)_Q = [A, B]_Q$.  
If necessary we extend the vector field 
$\dot\Gamma_h(t)$ along $\Gamma_h(t)$ 
as a genuine vector field around 
$\Gamma_h(t_0)$, this applies to 
the vector fields $A=X(\theta(t_0))$,  
$B=\dot\Gamma_h(t)-X(\theta(t_0))$, and 
the point $Q=\Gamma_h(t_0)$. 
Therefore we obtain 
$$
\nabla^P_{X(\theta(t_0))_{\Gamma_h(t_0)}}\dot\Gamma_h(t)
= 
\nabla_{X(\theta(t_0))_{\Gamma_h(t_0)}}\dot\Gamma_h(t)
= 
\left.\nabla_{\dot\Gamma_h(t)}
\dot\Gamma_h(t)\right\vert_{t=t_0}.  
$$ 

From the presentation 
$\dot\Gamma_h(t)=X(\theta_h(t))_{\Gamma_h(t)}$ 
and the above computations,  
we see 
\begin{eqnarray*}
\dot\theta_h(t_0)
\!&&\!\!\!\!\!\!
=
\left.\frac{d}{dt}\right\vert_{t=t_0}
\left\langle
\dot\Gamma_h(t), Y(\theta_h(t_0))_{\Gamma_h(t)}
\right\rangle
\\
=\!&&\!\!\!\!\!\!
\left\langle
\nabla_{\dot\Gamma_h(t_0)}\dot\Gamma_h(t),
 Y(\theta_h(t_0))_{\Gamma_h(t_0)}
\right\rangle
+
\left\langle
\dot\Gamma_h(t_0), 
\nabla_{\dot\Gamma_h(t_0)}Y(\theta_h(t_0))_{\Gamma_h(t)}
\right\rangle
\\
=\!&&\!\!\!\!\!\!
\kappa_g(t_0) 
+ 
\left\langle
X(\theta_h(t_0)), 
\nabla_{X(\theta_h(t_0))}Y(\theta_h(t_0))
\right\rangle
=
\kappa_g(t_0) .  
\end{eqnarray*}
\QED Propostion 4.13. 
\bigskip
\\
\proofof{Theorem 4.9}
We are ready to compute the linear holonomy of $W$ 
inside $\E/\W$.  
We start the computation without assuming that $\kappa$ is constant.  
With respect to the framing $\Theta $ and $\tilde Y$ of $\E/\W$, 
the infinitesimal action of $W$ is computed as 
\vspace{-3pt} 
$$
[W, \Theta] 
= 
[\tilde X + \tilde Z -(\kappa+1)\Theta,\,\Theta]
=
-\tilde Y \vspace{-3pt} 
$$
and \vspace{-3pt} 
\begin{eqnarray*}
[W, \tilde Y]
\!&=&\! 
[\tilde X + \tilde Z -(\kappa+1)\Theta,\,\tilde Y]
=\kappa(\tilde Z \tilde X) + 
(\tilde Y\cdot\kappa)\Theta
\\
\!&\equiv&\! 
\kappa(\kappa +1)\Theta 
+ 
(\tilde Y\cdot\kappa)\Theta
=
\{\kappa(\kappa +1)+(\tilde Y\cdot\kappa)\}\Theta
\quad (\mathrm{mod}\,\W). 
\end{eqnarray*}
This implies,  modulo $\W$ 
$$
\left.
\frac{d}{dt}\right\vert_{t=0}\exp(tW)_*(\Theta, \, \tilde Y) 
=
(\Theta, \, \tilde Y) 
\begin{pmatrix}
0 & -\{\kappa(\kappa +1)+(\tilde Y\cdot\kappa)\}
\\
1 & 0
\end{pmatrix}
. 
$$

Finally we assume hereafter 
the curvature $\kappa$ to be constant.   
Then the above matrix  reduces to   
$\ds 
\begin{pmatrix}
0 & -\kappa(\kappa +1)
\\
1 & 0
\end{pmatrix}
$, 
from which 
it is almost clear that the theorem holds.  
However we take a slightly closer look at what happens 
in each cases. 
\bigskip

\noindent
{\bf Parabolic case ``$\kappa=0,-1$'' :} 
\quad
If $\kappa\equiv 0$ or $-1$, we have
$\ds 
\exp(tW)_* 
=
\begin{pmatrix}
1 & 0
\\
t & 1
\end{pmatrix}
$.  
Therefore in $\E/\W$ along $\W$, 
$\tilde Y$ is invariant, without expanding nor contracting, 
and $\Theta\equiv\D/\W$ inclines to $\tilde Y$.  
There is no other invariant line field than 
$\langle\tilde Y\rangle$.  
They are of genuine-parabolic type. 
\bigskip

\noindent
{\bf Elliptic case ``$\kappa<-1$'' or ``$\kappa>0$'' :} 
\quad
If $\kappa<-1$ or $\kappa>0$, as $\kappa(\kappa +1)>0$, 
we take $K\Theta$ and $\tilde Y$ 
as a global framing 
in place of $\Theta$ and $\tilde Y$, 
where $K=\sqrt{\kappa(\kappa +1)}$.    
Then the matrix for the holonomy becomes  
$\ds 
\exp(tW)_* 
=
\begin{pmatrix}
\cos Kt & -\sin Kt
\\
\sin Kt & \cos Kt 
\end{pmatrix}
$,    
from which we see clearly that the system is of elliptic type. 
\bigskip

\noindent
{\bf Hyperbolic case ``$-1<\kappa<0$'' :} 
\quad
If $-1<\kappa<0$,  as $\kappa(\kappa +1)<0$, 
we take $K=\sqrt{-\kappa(\kappa +1)}$ and employ 
the same change of global framing. 
Then
Then the matrix for the holonomy becomes  
$\ds 
\exp(tW)_* 
=
\begin{pmatrix}
\cosh Kt & \sinh Kt
\\
\sinh Kt & \cosh Kt 
\end{pmatrix}
$.     
There fore we have two invariant subline bundles,  
$\ds l^u =
\left\langle 
K\Theta+\tilde Y 
\right\rangle
$ which is expanding,  
and 
$\ds l^s =
\left\langle 
K\Theta-
\tilde Y 
\right\rangle
$ which is contracting.  

$\D/\W=\langle \Theta \rangle$ is 
expanding and coming closer to $l^u$ 
as $t\to\infty$ while 
expanding and coming closer to $l^s$ 
as $t\to -\infty$ without passing through 
$l^u$ nor $l^s$.  
The system is of genuine-hyperbolic type.  
\smallskip
\\
\QED 
\bigskip 
Theorem 4.9. 
\Rem{
On the magnetic extension of surfaces 
of constant curvature $\kappa$, 
considering the symmetry, 
we see that the null-geodesics should project 
down to the surface to be curves with constant curvature.  
Let us look at the case of negative curvature. 
In order to compare the cases of different curvatures, 
we take homotethic transofromations 
(conformal transformation by multiplying constant)
so that 
all surfaces are of constant curvature $-1$   
and then look at projected curves. 
If the original surface is of constant curvarure 
$\kappa$, then multiplying $(-\kappa)^{-1/2}$ to the metric 
($(-\kappa)^{-1}$ to the Riemannian metric tensor $dh$ ) 
provides a surface of constant curvature $-1$. 
Then curves of constant geodesic curvature $-\kappa$ is 
transformed into those of $(-\kappa)^{1/2}$. 

$\kappa\equiv -1$ is the critical case. The curves of 
constant geodesic curvarture $1$ is nothing but 
the horocycles. 
If $\kappa<-1$, then the correponding curves on 
the universal covering $\Hyp$ after 
the conformal transformation are those of 
constant geodesic curvature $(-\kappa)^{1/2}>1$, 
which are closed circles, 
while in the case of  $-1<\kappa<0$ 
they are of geodesic curvature $ 0<(-\kappa)^{1/2}<1$ 
and are not compact. They have an intermediate character 
between  geodesics and horocycles.  

In the spherical case $\kappa>0$, 
always those curves are small circles on $S^2$, 
\EG for $\kappa\equiv 1$ then 
a circle which passes through the north pole and touches 
the equator is one of those. 

These observation explains Theorem 4.9 from a slightly different 
point of view. 
\bigskip
}

\subsection{Lorentz prolongation-III : Lorentzian surfaces${}^*$}

The product Lorentzian extension of a hyperbolic surface 
(or a surface of negative curvature in general) 
gives rise to one of the bi-Engel structure 
corresponding to the negative one 
of the bi-contact structure 
associated with the geodesic Anosov flow. 
The seek for a natural construction of Engel structure 
of the partner of the above mentioned one 
motivated the examples presented in this subsection.

Starting with 
a Lorentzian surface, 
we can perform the magnetic extension 
in order to obtain a Lorentzian 3-manifold and thus 
its Lorentz prolongation, 
while what we obtain so far still seems a bit mysterious. 
It is an interesting problem to analyze them.  
This is left to readers as a problem.

In this subsection we take 
the product of a Lorentzian surface 
$(\Sigma, dh)$ 
with $(S^1, -d\theta^2)$, so that we obtain 
a Lorentzian 3-manifold $(V, dg)$. 
We do not assume 
$\Sigma$ to be compact, 
while for the sake of simplicity we assume 
the orientations not only of $\Sigma$ but also 
of the positive and negative directions. 
It means the following. 
At each point $\sigma\in\Sigma$ 
the set 
$S^c=S^c (T_\sigma\Sigma) 
= \{u\in T_\sigma\Sigma 
\vert dg(u, u)=c\}$ is 
a hyperbola if $c\ne 0$ or a pair 
of straight lines crossing at the origin if $C=0$. 
The orientation of the positive or negative directions 
implies the choice of a connected component 
$S^{\pm1}_0$ of $S^{\pm 1}$.  
We also assume that 
the surface is oriented and 
the pair 
$\langle u_+, u_- \rangle$ 
($u_\pm \in S^{\pm 1}_0$) 
forms an oriented basis. 

In this situation 
the structural group $\mathit{O}(1;1)$ 
reduces to the connected component 
$\mathit{SO}(1;1)_0$ 
of the identity and 
the set of oriented Lorentzian orthonormal 
frame 
$\ds \cup_{\sigma\in\Sigma} 
\{\langle u_+, u_- \rangle\, \vert\, 
u_\pm \in S^{\pm 1}_0\} 
$ 
is naturally identified with 
the associated principal 
$\mathit{SO}(1;1)_0$-bundle over $\Sigma$.  
If we take an oriented Lorentzian orthonormal basis 
$u_\pm$ at a point $\sigma\in\Sigma$ as above, 
then the other (negative) unit vectors in $S^{\pm 1}_0$  
are indicated as $\cosh t u_+ + \sinh t u_- $ 
or as $\sinh t u_+ + \cosh t u_- $.

Like in the Riemannian case, using the Levi-Civita 
connection, we have the tautological horizontal 
vector fields $X$ and $Y$ corresponding to 
 $u_+$ and $u_-$ on the unit tangent bundle 
$S^1(T\Sigma)$. 
The vertical vector field which generates 
the action of $\R=\{t\}$ in the above sense 
is denoted by $Z=\frac{\partial}{\partial t}$. 
Then we have the structural equations  
$$
[Z,X]=Y\,,\quad 
[Z,Y]=X\,,\quad 
[X,Y]=\kappa Z\
$$ 
where $\kappa$ denotes the Lorentzian 
curvature of $(\Sigma, dh)$. 
The positive geodesic flow is generated by $X$, 
the negative one by $Y$. 
\smallskip

The three dimensional Lorentzian metric 
$dg$ has signature $(1,2)$, namely, 
one positive dimension and two negative dimensions. 
At each point $v\in V$ the set $NC(T_vV)$ of 
null-lines in the tangent space $T_vV$ is again a circle. 
On the other hand, in this construction 
the circle contains 
two special points 
$\ell_+$ 
and 
$\ell_-$ 
which are the fixed points of the involution  
$\iota : v=(\sigma, \theta) \mapsto (\sigma, -\theta)$.  
These two lines are those who already exist in 
$T_\sigma\Sigma$ as null-lines 
$\ell_+=
\langle 
u_+ + u_-
\rangle  
$ 
and 
$\ell_-=
\langle 
u_+ - u_-
\rangle  
$. 
The complement of these two points consists of 
two open arcs which are exchanged by the involution. 
So we take the only one open arc in this construction. 
This open arc $NC^+=NC^+(T_\sigma\Sigma)$ 
is also regarded as the real line 
$\R=\{t\}$ by the correspondence 
$
t \leftrightarrow
\langle 
\cosh t u_+ + \sinh t u_- + \Theta
\rangle
$ 
after fixing an oriented Lorentzian ON basis 
$u_\pm$.   
Here again $\Theta$ 
denotes $\ds \frac{\partial}{\partial \theta}$.   
Therefore like in 4.1, our 4-manifold 
$M=NC^+(TV)$ is naturally identified with 
$S^1(T\Sigma)\times S^1$ by 
$((\sigma, \theta), u_++\Theta)$ 
$\leftrightarrow$ 
$((\sigma, u_+),\theta)$.  
The natural lift of the null-geodesics are 
generated by $W=X+\Theta$ and the Engel structure 
is given by 
$\D=\langle W,\, Z  \rangle$.  
and we obtain the even contact structure 
$\E=\langle W,\, Z, \, Y  \rangle$.  
As $\Theta$ commutes with $X$, $Y$, and $Z$, 
this is also verified by the commutation relation. 

Now we focus our attentions to a specific model. 

\begin{ex}{\rm
{\bf(de Sitter space)}
\quad 
Let us take the de Sitter space 
$dS_2=\{
(s_1,s_2,s_3)\in \R^3\,\vert\,
s_1^2+s_2^2-s_3^2=1
\}$ as $\Sigma$.   
which is also identified with 
$\mathit{PSL}(2;\R)/\{\exp{tl}\,\vert\,
t\in\R\}
$. 
Here we are following the notations in 
Example 4.4.  
More precisely, 
$\R^3$ and its standard basis 
correspond to 
the Lie algebra $\mathit{psl}(2,\R)$ 
and their basis $h$, $l$, $k$ 
and the $(2,1)$-type metric on $\R^3$ 
to the adjoint invariant metric 
(\IE the Killing form) 
$\mathrm{tr}(\mathit{ad}(\cdot)\circ\mathit{ad}(\cdot))
=2\mathrm{tr}(\cdot\times\cdot)$. 
Then as is well-known, the principal 
$\mathrm{SO}(1,1)_0$-bundle 
which coincides with the unit tangent bundle 
$S^1(TdS_2)$ over $dS_2$ is identified with 
the principal $\R$-bundle 
$\mathit{PSL}(2;\R)\to 
\mathit{PSL}(2;\R)/\{\exp{tl}\,\vert\,
t\in\R\}
=dS_2$.  
On the total space 
$S^1(TdS_2)=\mathit{PSL}(2;\R)$ 
the canonical vector fields $X$, $Y$, and $Z$ are now 
nothing but the left invariant vector fields 
$h$, $k$, and $l$. From the commutation relation, 
we see that $\kappa\equiv -1$. 

If we take a closed hyperbolic surface 
$\Gamma\setminus\Hyp^2$ with 
$\pi_1(\Gamma\setminus\Hyp^2)\cong\Gamma 
\subset\mathit{PSL}(2;\R)$, 
$\Gamma$ acts as orientation preserving isometry on 
$\Gamma\setminus\Hyp^2$, the whole construction 
is invariant with respect to this action, 
namely, to the left translation of $\Gamma$ 
to $\mathit{PSL}(2;\R)\times S^1$. 
The action on the second factor is trivial.  
Thus we obtained the 4-manifold 
$M'=S^1(T\Gamma\setminus\Hyp^2)\times S^1$ 
and an Engel structure 
$\D'=\langle h+\Theta, l \rangle$, 
the Cauchy characteristic 
$\W'=\langle W=h+\Theta \rangle$, 
and the even contact structure 
$\E'= \langle h+\Theta, l, k \rangle$.   

As mentioned in Remark 4.3, $M'$, $\W'$, and $\E'$ are 
exactly the same as in Subsection 4.1 for 
a hyperbolic surface $\Gamma\setminus\Hyp^2$ 
and the Engel structure obtained here is 
one of the bi-Engel structure (\cite{KV}) 
which corresponds to 
the contact structure $\xi_+$ explained in Remark 4.3.  
}
\end{ex}

\Prob{1)\quad 
We do not know when we start with the magnetic 
extension of a Lorentzian surface and adopt 
the above construction 
with $NC^+$ whether if there exists some good 
discrete group action which yields a compact quotient. 
\\
2)\quad
If we could have 
a nice compactification in the above problem, 
it should be interesting to look at the dynamics and 
look for the parabolic ones.  
\\
3)\quad
Also, even in the above example, we have not yet 
understood  
the meaning of or the geometry corresponding to 
the compactification 
of the subspace $NC^+$ of null-lines 
by $\{\ell_+, \ell_-\}$.   

Looking both on the product extension and 
the magnetic extension concerning these problems 
might be of some interest. 
}

\subsection{Pre-quantum prolongation${}^*$}
 Our aim in this subsection is to give examples of 
pre-quantum prolongations and look for 
more Engel structures of 
hyperbolic or possibly parabolic type.   
Basically we follow the notations in Subsection 1.4. 
\begin{ex}[Nil-Solv hybrid]
{\rm\quad 
Consider a solvable 3-manifold 
$V^3$ which fibers over the circle $S^1$ 
with fibre $T^2$ and a hyperbolic monodromy 
$\varphi\in\mathit{SL(2;\R)}$ ($\mathrm{tr}\,\varphi >2$). 
The standard area form of area  $1$ is invariant under the 
monodromy $\varphi$, we can take 
a closed 2-form $\omega$ whose restriction to any fibre 
is the area form, which represents  
the cohomology class $\alpha\in H^2(V;\Z)\cong \Z$ 
corresponding to $1\in\Z$.  
We take the the volume form 
$\vol=d\theta\wedge\alpha/2\pi$.  
Then the suspension flow $\ul{W}$ is a lift of 
$\ds 2\pi\frac{d}{d\theta}$ and is the Poincar\'e dual 
to $\alpha$ with respect to $\vol$, \IE 
$\alpha=\iota_{\ul{W}}\vol$. 
$\ul{W}$ is one of the standard suspension Anosov flows. 
Now we take the 
bi-contact structure $\xi_\pm$ associated with this Anosov flow 
(\cite{KV}, \cite{Mi}). 
$\ul{W}$ is Legendrian for both of $\xi_\pm$.  
Performing the pre-quantum prolongation process, 
we obtain a bi-Engel structure $\D_\pm$ 
whose even contact structure is the horizontal 
space defined by the connection 1-form whose 
curvature form is exactly $\omega$. 
The resulting 4-manifold $W$ is the circle bundle 
over $V$, whose restriction to each fibre $T^2 \subset V$ 
is a nilpotent 3-manifold $\mathit{Nil}^3(1)$ of euler class $1$. 

The vector field $W$ which spans 
the Cauchy characteristic $\W$ 
is a lift of $\ul W$ and preserves 
the length of the fibre circle. 
Therefore the dynamics is nothing 
almost the same as that of the Anosov flow 
$W$.    

\smallskip

This is one way to construct a bi-Engel structure 
introduced in \cite{KV}. 
We can also start from the 4-dimensional Lie group $G$ 
which is the central extension of the 
3-dimensional solvable Lie group $\mathit{Solv}^3$. 
Namely, $G$ is obtained by 
$$
0\to\R^2 \to \mathit{Solv}^3 \to \R \to 0\,, \qquad 
0\to\R \to G \to \mathit{Solv}^3 \to 0
$$ 
where the first exact sequence is a semi-direct product 
by the action of $\varphi^t\in \R$ on $\R^2$ and the
second one is the unique (up to constant) 
non-trivial central extension.  Then we can find 
a co-compact lattice in $G$.  

This Lie group is also presented as a non-central extension 
of the 3-dimensional Heisenberg Lie group $H_\R$ 
by an automorphism $\tilde\varphi$ which preserves 
the integral lattice  $H_\Z$ where 
$$
H_\R=\left\{
 \begin{pmatrix}
 1&x&z\\0&1&y\\0&0&1\\
\end{pmatrix}
\left\vert\right. x,y,z\in \R
\right\}
\,\,\,
\supset
\,\,\,
H_\Z=\left\{
\,\,*\,\,
\left\vert\right. x,y,z\in \Z
\right\}
$$
and $\tilde\varphi$ acts as 
$\begin{pmatrix}
x\\y\\
\end{pmatrix}
\mapsto
\varphi
\begin{pmatrix}
x\\y\\
\end{pmatrix}
$. 
}
\end{ex}
\begin{ex}[Geodesic flow of Riemannian surface]
{\rm\quad
Let $(\Sigma, dg)$ be a Riemannian surface 
with no-where zero 
curvature $\kappa$. 
Consider the unit tangent bundle $p :V=S^1(T\Sigma)\to \Sigma$, 
the geodesic flow on $V$ generated by the horizontal 
canonical vector field $X$, and the Levi-Civita connection 
whose connection 1-form is denoted by $\beta$.  
$\xi=\mathrm{ker}\,\beta$ is the horizontal plane field 
to which $X$ is tangent. 
If the curvature $\kappa$ is everywhere 
positive \RESP{negative} 
$\xi$ is a negative \RESP{positive} contact structure, 
because $d\beta=-\kappa p^*d\mathrm{area}_\Sigma$.   
Because $X$ is the Reeb vector field 
of the canonical contact 1-form 
(\IE of the associated Liouville 1-form) $\lambda$, 
$X$ is an exact vector field in the sense of asymptotic cycle. 
We adopt $X$ as $W$ for the pre-quantum prolongation.  
The invariant volume is $\lambda\wedge d\lambda$ and 
the closed 2-form to which the pre-quantization is performed is 
$\omega=d\lambda$, 
so that the corresponding cohomology class $\alpha$ 
can be any of the torsion part 
$\mathrm{Tor} (H^2(V;\Z))$. 

The resulting manifold $M$ is the circle bundle over $V$ 
whose euler class is $\alpha$. 
The action of the Cauchy characteristic $W$ on $\E/\W$ 
is almost same as the dynamics of the geodesic flow on $V$. 
Therefore if the curvature $\kappa$ is everywhere positive, 
the Engel structure is of elliptic type, and if $\kappa$ is 
everywhere negative, it is of hyperbolic type. 
Especially  if we take the cohomology class 
$\alpha=0=[d\lambda]$, 
on the resulting 4-manifold $M=V\times S^1$, 
with respect to this product structure we can take 
the connection form to be $d\theta+\pi^*\lambda$.  
Therefore the Cauchy characteristic $\W$ is generated by 
the vector field 
$\ds W=X-\frac{\partial}{\partial \theta}$.  
Therefore the Engel structure we obtained coincides 
with one in Proposition 4.1 which is obtained by 
the suspension of the time $2\pi$ map of the geodesic flow. 
\medskip

For a closed oriented surface 
$\Sigma=\Sigma_g$ of genus  $g\ne 1$, 
we have non-unique candidates for the resulting manifold $W$ 
on which the Engel structure 
because $\mathrm{Tor}(H^2(V;\Z))\cong\Z/$
{\small $\!\!\vert2g-2\vert$}.  

For example, if we assume that $\kappa>0$, 
$\Sigma$ is topologically $S^2$, $V$ is diffeomorphic to 
$\R P^3$, and $H^2(V;\Z)\cong\Z/2$. 
If we take $\alpha=0$, $M\cong V\times S^1$ and 
$\pi_1(M)\cong \Z\times \Z/2$.  
On the other hand, if we take $\alpha\ne0\in\Z/2$, 
the fiber bundle structure $S^1\to M \to \R P^3$ 
and the euler class imply that 
$\pi_1(M)\cong\Z$ 
(the $\pi_1$'s of the fibration is an extension : 
$0\to\Z\overset{\times 2}{\to} \Z \to \Z/2\to 0$). 
In fact $M$ is diffeomorphic to $U(2)$ and 
even to $S^3 \times S^1$.  

Unfortunately, even the manifold is different, 
the Engel structure is not so new. 
The natural projection 
$M\to S^1$ to the circle of the half length 
$S^1/${\small $\Z/2$} tells that $M$ is the mapping 
torus of the antipodal map $\tau$ of $S^3$, 
which is isotopic to the identity of $S^3$.  
Therefore $M$ is diffeomorphic to the product. 
The antipodal map $\tau$ is the deck transformation 
of the universal covering $S^3 \to \R P^3$ which pulls 
the geodesic flow $\phi_t$ 
back to $S^3$ as the flow 
$\tilde\phi_t$ 
and the Liouville contact structure 
$\xi_0=\mathrm{ker}\lambda$ to $\tilde\xi_0$ as well. 
Therefore $\tau$ and $\tilde\phi_t$ commute 
to each other. 
Our Engel structure is also obtained 
by the suspension construction given in Subsection 1.5 
with respect to $V=S^3$, $\tilde\xi_0$, and 
$\varphi=\tau\circ\tilde\phi_{2\pi}$.

Like in the case of $\kappa>0$, for a surface of 
$g>1$ and $\kappa<0$, 
consider finite converings in the fibre direction 
of $V \to \Sigma_g$, we see that 
more or less similar situations appear.    
}
\end{ex}
\begin{ex}[Propellor constructions-I]
{\rm\quad
This class of examples are a generalization of 
the Nil-Solv hybrid example. 
For the suspension construction, we need 
a contact structure on 3-manifold which gives rise to 
the even contact structure, 
while in the pre-quantum 
prolongation, the non-vanishing closed 2-form 
replaces the roll of the non-integrability 
which assures the bracket generation $[\E,\E]=TM$. 

Let $V$ be the mapping cylinder 
of a linear automorphism $\varphi\in \mathit{SL}(2;\R)$ 
of $T^2$.  
It fibers over the circle and 
is also presented as the 
quotient of $\tilde V=T^2\times S^1\ni (x,y,t)$ 
($x, y, \in \R/\Z$) by the identification 
$(x,y, t+1)\sim(x', y',t)$ 
with 
$
\begin{pmatrix}x'\\y'\\
\end{pmatrix}=
\varphi\begin{pmatrix}x\\y\\
\end{pmatrix}
$.  
$\frac{\partial}{\partial t}$ 
on $\tilde V$ descends to $\ul W$ on $V$.  

The propellor construction of contact structure 
is to choose a linear foliation 
$\ds \tilde l_t=\left\langle a(t)\frac{\partial}{\partial x}
+
b(t)\frac{\partial}{\partial y}\right\rangle$  
on each $T^2\times \{t\}$ 
in such a way that $\l_t$ rotates in 
positive or negative direction when $t\in\R$ increases 
and satisfies 
$l_{t+1}=\varphi_*l_t$.  
Then $\ds \tilde\xi=l\oplus
\left\langle
\frac{\partial}{\partial t}
\right\rangle
$
is a contact structure on $\tilde V$ which descends 
to $\xi$ on $V$.

The fiberwise area form $\omega=dx\wedge dy$, 
the standard volume form $\vol = dx\wedge dy\wedge dt$ 
are also well-defined on $V$ and satisfies 
$\omega=\iota_{\ul W}\vol$ which is a closed 2-form. 
The cohomology class $\alpha=[\omega]$ is 
the Poincar\'e dual to the 
suspension circle and is integral. 
Therefore the pre-quantum prolongation can be 
performed to obtain an Engel structure on $M$ which fibers 
over $V$ with euler class $\alpha$. 

According to $\varphi$ being elliptic (including the case 
$\varphi$ is the identity), parabolic, or hyperbolic, 
the Engel structure is of elliptic, 
(genuine or trans-) parabolic, 
or (genuine or trans-) hyperbolic type.   
}
\end{ex}
\begin{ex}[Propellor constructions-II]
{\rm\quad
In the above propellor construction, 
we choose a non-singular vector field 
$\ds \widetilde{\ul W}= a(t)\frac{\partial}{\partial x}
+
b(t)\frac{\partial}{\partial y}$  
in such a way that not only $\tilde l$ is invariant under 
the monodromy but also the vector field itself is invariant, 
namely it satisfies 
$\ds \varphi\begin{pmatrix}
a(t)\\b(t)\\
\end{pmatrix}
=
\begin{pmatrix}
a(t+1)\\b(t+1)\\
\end{pmatrix}
$, 
so that $\widetilde{\ul W}$ descends to $\ul W$ on $V$. 
$\ul W$ is Legendrian and preserves $\vol$.  

In hyperbolic case or in elliptic case except for 
the identity, 
because through the projection 
$H_1(V;\Z)\cong H_1(S^1;\Z)\cong\Z$, 
$\omega=\iota_\vol$ is always exact.  
In parabolic case, if the eigenvalue of $\varphi$ 
is equal to $-1$, then again we have
$H_1(V;\Z)\cong H_1(S^1;\Z)\cong\Z$ and $\alpha=0$,  
while if it is equal to $1$, 
we have rank $1$ choice for $\alpha$, which depends 
on the choice of $\ul W$.  
In the case where $\varphi$ is the identity, 
we have choice of rank $2$. 
Multiplying appropriate non-vanishing smooth function 
on $S^1$ to $\ul W$ achieves the choice. 

It is not hard to see that in any case the resulting 
Engel structure is of parabolic type, 
while it is not surprizing because 
on the cyclic covering of $M$ which covers 
the cyclic covering $\tilde V \to V$, 
all of them look alike. 
}
\end{ex}

\section{Problems and discussions}

To close this note we collect problems concerning 
the topics discussed in this note.  
Some of them have already been mentioned before 
and some other may be accompanied with short discussions. 
Also some are concrete and others are rather vague. 
\bigskip

We have not discussed the symmetry of Engel manifolds.  
One reason is that for an Engel structure 
on a closed 4-manifold $M$, in general 
the symmetry is expected to be very small. 
First of all, if a diffeomorphism $\Phi:M\to M$ 
preserves the Engel structure $\D$, necessarily it also preserves 
$\W$. 
In generic cases it implies $\Phi$ preserves each orbit of $\W$. 
Then inside each $\W$-orbit the position of $\D/\W$ in $\E/\W$ 
almost determines the point of the orbit. 

The Cartan prolongation is an exceptional case in this sense. 
If we start from a contact 3-manifold $(V,\xi)$ and obtain 
its Cartan prolongation $(M,\D)$, the symmetry of 
 $(M,\D)$ is exactly that of  $(V,\xi)$, \IE 
naturally we have 
$\mathit{Diff}^\infty(M,\D)\cong\mathit{Diff}^\infty(V,\xi)$.  

If a closed Lorentzain 3-manifold $(V, dg)$ admits an isometry, 
it automatically lifts to a symmetry of its Lorentz prolongation. 

There are  many other closed Engel manifolds 
for which we have continuous symmetries. 

\Prob{ 
Classify all the closed Engel manifolds with continuous symmetry. 

In particular, study the case where 
in the support of continuous symmetry 
the $\W$-orbits are not closed.  
}

On the other hand, the 2-jet space $J^2(1,1)$ 
and the standard Engel structure $\D_0$ 
with the standard Engel-Darboux coordinates   
admits a big symmetry group, which is naturally 
isomorphic to the group of contactomorphisms 
of the standard contact structure $(\R^3, \xi_0)$. 
Here we follow the notations in [EF] in Subsection 1.1.  
In side this group we can find a 4-dimensional 3-step 
nilpotent Lie group $G_E$ which acts transitively and freely 
on  $\R^4=J^2(1,1)$. 

The vector fields 
$W$, $X$, $Z$, and $Y$ on $\R^4$ generate a 4-dimensional 
nilpotent Lie algebra $\g_E$ and thus the corresponding 
nilpotent Lie group $G_E$. 
In fact they satisfy the commutation relations 
$[W, X]=Z$, $[Z,X]=Y$, and others are trivial. 
These relations exactly corresponde to the flag generation 
$\W\subset\D\subset\E\subset T\R^4$. 
However, apparently among these vector fields 
$W$ and $Z$ are out of symmetries of $\D_0$. 
If we identify $\R^4=J^2(1,1)$ with the Lie group $G_E$ 
and $\g_E$ is the set of left invariant vector fields on $G_E$, 
the full flag $\W\subset\D\subset\E$ is understood as 
left invariant fields on $G_E$. In this formulation 
the action, \EG $\exp(tW)$ generated by $W$ is a right 
translation on $G_E$ and thus necessarily preserve 
left invariant fields.  
Relying on the left action which preserves the Engel structure 
the element $\Phi_{dcba}\in G_E$ is indicated as 
$\Phi_{dcba}:(x,y,z,w)\mapsto 
(x+d,\, y+\frac{1}{2}ax^2+bx+c,\, z+ax+b,\, w+a)$  
where $\Phi_{dcba}(0,0,0,0)=(d,c,b,a)$ for 
$(d,c,b,a)\in\R^4$. 

The association of the Lie algebra $\G_E$ or the Lie group $G_E$ 
to an Engel structure seems to depends on the choice of 
local Engel-Darboux chart. However, even at any point and its 
neighbourhood, the Lie algebra structure seems to survive as 
the generation of the flag. 
Therefore $\g_E$ msut be very fundamental.  

\Prob{An Engel structure may have continuous symmetry like 
in the case of Lorentz prolongation from 
a Lorentzian 3-manifold with isometry group of positive dimension. 
Formulate relations between such global symmetry to 
$\g_E$ or if possible to $G_E$. 
}

\Prob{In which sense dess the {\it h}-principle for Engel structures 
hold? Describe the path-components of the set of Engel structures. 
}

\begin{prob}{\rm
( See 2) in Remark 1.7, \cite{Ch}, and \cite{SY}.)
\quad
Formulate the W\"unschmann invariant 
for a line field ${\mathcal L} \subset \D$ 
to be the fibre of the Lorentz prolongation in a local sense 
in a manner which is suitable for the following purpose.  
For a given Engel structure 
study the existence or an obstruction to the existence of 
a line field ${\mathcal L} \subset \D$ 
with vanishing W\"unschmann invariant, 
or the deformability to such one in a 1-parameter family of 
Engel structures of a given one.  
}
\end{prob}

\Prob{
Which kind of informations on Lorentzian manifold 
can we deduce from their Lorentz prolongation?
How about the same question for the above generalized case, 
Engel structures equipped with line fields 
${\mathcal L} \subset \D$ 
with vanishing W\"unschmann invariant?
}

\Prob{Inaba's accessible set 
should be reconsidered 
from the coordinate-free point of view.  

It might be possible to distribte 
the germ of the accessible set along a $\W$-curve 
of each point 
as a kind of field on the Engel manifold.  
Formulate this notion and relate it to the 
global and dynamical structure of Engel 
manifolds. 
}

Consider the Minkowski space $\R^{2,1}=\{(r,s,t)\}$ 
with $dg=dr^2 + ds^2 -dt^2$ and its Lorentz prolongation. 
We can arrange an Engel-Darboux coordinate 
in such a way that $
r=\frac{1}{\sqrt 2}(w-y),\, s=z,\,
t=\frac{1}{\sqrt 2}(w+y)$, 
and the positive ineteriot $A_+$ of Inaba's accessible set 
from the origin coincides with the interior of the causal set of 
the origin in $\R^{2,1}$.  In this case, 
the fibre direction of the Lorentz prolongation is exactly 
$\langle X \rangle$ in [EF]. 
The function $z^2-2wy=r^2+s^2-t^2$ is a 1st integral of 
the vector field $X$.   

\Prob{Study the relationship between this fact and the problems 
5.3 - 5.5 above. 
}

\begin{prob}
{\rm
(cf. Subsection 3.4.) \quad 
Give a proof of Theorem \ref{null-geodesic}  
relying on the causality and the Bryant-Hsu rigidity.  
}
\end{prob}

\Prob{
Develop the study of (the group of) contactomorphisms 
of 3-dimensional contact structures 
in order not only for the construction 
but also for classification problem of Engel structures. 
}

\Prob{
Study the Engel structure of Lorentz prolongation from 
one more family of Lorentzian extension of 
Riemannian surface;  
$(\Sigma, dh)\times (S^1, -f(\theta)(d\theta)^2)$ 
for a positive non-constant function $f$.  
$S^1$ can be repalced with any 1-manifold. 
From the point of view of Engel structures, 
it is equivalent to study 
$(\Sigma, f(\theta)^{-1}dh)\times (S^1, -(d\theta)^2)$.  
}

\Prob{In the construction of compact Lorentzian 3-manifold
from product extension 
of Lorentzian surfaces in Subsection 4.3, 
we did not take the whole of the null circle bundle 
but a portion of it 
corresponding to an open interval in null circle. 
After taking a quotient by a discrete group action, 
the boundary disappeared. 
Which kind of informations do the boundary or 
the compactification carry? 
}

\Prob{
Study the 
magnetic extensions of Lorentzian surfaces 
and associated Engel structures. 
}

\Prob{
By definition, paraboic/hyperbolic Engel structures  
require only two dimensional foliations to which 
the Cauchy characteristic $\W$ is tangent, 
while in any example in this note it is raised to 3-dimensional 
ones! 
Does there exist examples which can not be raised, or is it 
a natural consequence? 

How about the Engel structures coming from Eliashberg-Thurston's example 
based on laminations \cite{ET}?
}

\begin{Prob}
{\rm $\!\!\!\!\!\!$(Discussions on Example 4.20 and Proposition 4.1.) 
\quad 
What is the limit of Example 4.20  with $\alpha=0\in H^2(V; \Z)$ 
when $\kappa$ tends to $0$? 
In the family of 
isomorphic examples in Proposition 4.1,  
$\kappa=0$ does not matter and we have an Engel structure. 
\medskip 

For a Riemannian surface,  in Subsection 4.1 
first we performed a (special) relativistic procedure 
(taking the product with $(S^1, -d\theta^2)$) 
and then we took a phase space 
for lightlike motions and obtained 
an Engel structure.   

On the other hand, in Example 4.20  with $\alpha=0\in H^2(V; \Z)$, 
we first took the 
phase space for motions of fixed kinetic energy 
on the surface and then we perform the pre-quantization, 
then again we obtained the same Engel structure. 

Does there any Physical significance of the 
coincidence of the results of these procedures?  
(Taking phase spaces conjugates 
the special relativistic procedure 
and quantum one. )
If we look into these examples, the coincidence is 
natural and nothing is mysterious...  

If there is some meaning, then does there exists 
further implication of the non-uniqueness 
of the case where $H^2(V; \Z)$ has torsions? 
\medskip

We have also the third construction for the same Engel structure, 
namely the suspension 
by the geodesic flow (at certain time) 
(the method of Subsection 1.5) of the Liouville 
contact structure.  
This construction sounds more purely mathematical, 
but in this case we need to put a plane field $\D$ somehow 
by hand.   
}
\end{Prob}

\Prob{Provide the study of Engel structures with singularities 
with interesting examples and a framework or guiding principle.  }

\begin{flushright}
Yoshihiko MITSUMATSU 
\vspace{3pt}
\\
{\small Department of Mathematics, Chuo University
\\
1-13-27 Kasuga Bunkyo-ku, Tokyo, 112-8551, Japan
\\
e-mail : yoshi@math.chuo-u.ac.jp
\vspace{3pt}
\\
}
\end{flushright}

\end{document}